\documentclass[11pt]{article}
\usepackage[T1]{fontenc}
\usepackage{natbib}
\usepackage{amsmath,amsfonts,amssymb,amsthm}
\usepackage{a4wide}
\usepackage{color}
\usepackage{epsfig,psfrag}

\newtheorem{dfn}{Definition}[section]
\newtheorem{thm}[dfn]{Theorem}
\newtheorem{lem}[dfn]{Lemma}
\newtheorem{cor}[dfn]{Corollary}
\newtheorem{rem}[dfn]{Remark}
\newtheorem{prop}[dfn]{Proposition}

\newenvironment{dem}{\vskip 2mm\noindent {\it Proof} :}
                    {\hfill $\square$ \vskip 2mm \noindent}
 
\DeclareMathOperator{\cyl}{cyl}
\DeclareMathOperator{\hyp}{hyp}
\DeclareMathOperator{\card}{card}

\DeclareMathOperator{\ad}{ad}

\DeclareMathOperator{\lleft}{left}
\DeclareMathOperator{\rright}{right}

\newcommand{\eps}{\varepsilon}

\def\PP{\mathbb{P}}
\def\RR{\mathbb{R}}
\def\EE{\mathbb{E}}
\def\NN{\mathbb{N}}
\def\ZZ{\mathbb{Z}}
\def\E{\mathcal{E}}
\def\D{\mathcal{D}}
\def\H{\mathcal{H}}
\def\N{\mathcal{N}}

\def\M{\mathcal{M}}

\def\wt{{\widetilde{\theta}}}
\def\F{\mathcal{F}}
\def\G{\mathcal{G}}

\def\LL{\mathbb{L}}
\def\va{\vec{v}(\theta)}
\def\vb{\vec{v}(\widetilde{\theta})}
\def\vc{\vec{v}^{\bot}(\theta)}
\def\vd{\vec{v}^{\bot}(\widetilde{\theta})}

\def\II{\mbox{ 1\hskip -.29em I}} 

\begin{document}

\thispagestyle{empty}


\title{Law of large numbers for the maximal flow through tilted cylinders
  in two-dimensional first passage percolation}

\date{}
\author{}
\maketitle

\begin{center}
\vskip-1cm {\Large Rapha\"el Rossignol\footnote{Rapha\"el Rossignol was supported by the Swiss National
    Science Foundation grants 200021-1036251/1 and 200020-112316/1.}}\\
{\it Universit\'e Paris Sud, Laboratoire de Math\'ematiques, b\^atiment 425,
91405 Orsay Cedex, France}\\
{\it E-mail:} raphael.rossignol@math.u-psud.fr\\
\vskip0.5cm and\\
\vskip0.5cm {\Large Marie Th\'eret}\\
{\it \'Ecole Normale Sup\'erieure, D\'epartement Math\'ematiques et
Applications, 45 rue d'Ulm, 75230 Paris Cedex 05, France}\\
{\it E-mail:} marie.theret@ens.fr
\end{center}



\begin{abstract}
Equip the edges of the lattice $\ZZ^2$ with i.i.d. random
capacities. We prove a
law of large numbers for the maximal flow crossing a rectangle
in $\RR^2$ when the side lengths of the rectangle go to infinity. The value
of the limit depends on the asymptotic behaviour of the ratio of the height
of the cylinder over the length of its basis. This law of large numbers
extends the law of large numbers obtained in \cite{GrimmettKesten84} for
rectangles of particular orientation.
\end{abstract}



\section{Introduction}

The model of maximal flow in a randomly porous medium with independent and
identically distributed capacities has been introduced by
\cite{Chayes} and  \cite{Kesten:flows}. The purpose of this model is to understand the behaviour of the maximum
amount of flow that can cross the medium from one part to another. 

All the precise definitions will be given in section \ref{chapitre6sec:notations}, but let
us draw the general picture in dimension $d$. The random medium is
represented by the lattice $\ZZ^d$. We see each edge as a microscopic pipe which the fluid
can flow through. To each edge $e$, we attach a non-negative capacity $t(e)$ which
represents the amount of fluid (or the amount of fluid per unit of time) that
can effectively go through the edge $e$. Capacities are then supposed to be
random, identically and independently distributed with common distribution
function $F$. Let $A$ be some hyperrectangle in $\RR^d$ and $n$ an
integer. The portion of medium that we will look at is a box $B_n$ of basis $nA$ and of height $2h(n)$, which
$nA$ splits into two boxes of equal volume. The boundary of $B_n$ is thus split
into two parts, $A_n^1$ and $A_n^2$. There are two protagonists in
this play, two types of flows through $B_n$: the maximal flow $\tau_n$ for which the fluid can enter the box
through $A_n^1$ and leave it through $A_n^2$, and the maximal flow $\phi_n$ for which the fluid enters $B_n$ only through its
bottom side and leaves it through its top side. The first quality of
$\tau_n$ is that it is (almost) a subadditive quantity, whereas
$\phi_n$ is not. The main question now is: ``How do $\phi_n$ and
$\tau_n$ behave when $n$ is large ?''.

In this paper, we shall understand this question as ``Is there a law
of large numbers for $\phi_n$
and $\tau_n$ ?'', and let us say that such results do indeed already
exist. However, it is
important to stress that the orientation of $A$ plays an
important role in these results. Indeed, the first ones
were obtained for ``straight'' boxes, i.e., when $A$ is of
the form $\prod_{i=1}^{d-1} [0,a_i] \times \{0\}$. Especially
concerning the study of $\phi_n$, this simplifies considerably the
task. Let us draw a precise state of the art. The law of large numbers for
$\tau_n$ were proved under mild hypotheses: in \cite{Kesten:flows} for
straight boxes and in \cite{RossignolTheret08b} for general boxes. These
results follow essentially from the subadditivity property already alluded to. Suppose that $t(e)$ has finite
expectation, $\vec{v}$ denotes a
unit vector orthogonal to a hyperrectangle $A$ containing the origin
  of the graph, and
$h(n)$ goes to infinity. Then there is a
function $\nu$ defined on $S^{d-1}$  such that:
$$ \nu(\vec{v}) \,=\, \lim_{n\rightarrow \infty} \frac{\tau(nA,
  h(n))}{\H^{d-1}(nA)} \qquad \textrm{a.s. and in }L^1\,, $$
where $\H^{d-1}(nA)$ is the $(d-1)$-dimensional Hausdorff measure of
$nA$. If the height function $h(n)$ is negligible compared to $n$, $\phi_n$
satisfies the same law of large numbers as $n$ (see for example
\cite{RossignolTheret08b}). Otherwise, the law of
large numbers for $\phi_n$ was proved only for straight boxes, with suboptimal
assumptions on the height $h$, the moments of $F$ and on $F(\{0\})$,
in \cite{Kesten:flows}. In dimension 2, this was first studied in \cite{GrimmettKesten84}. The assumption on $F(\{0\})$ was optimized in \cite{Zhang} and
\cite{Zhang07}. The assumptions on the moments of $F$ and the height $h$
have been improved in \cite{RossignolTheret08b}. A specificity of the
lattice $\ZZ^d$, namely its invariance under reflexions with respect
to integer coordinate hyperplanes, implies that the law of large
numbers is the same for $\phi_n$ and $\tau_n$ in straight cylinders
(provided $\log h(n)$ does not grow too fast).

Summarizing, $\tau_n$ is fairly well studied concerning laws of large
numbers, but for $\phi_n$, nothing is known when the boxes are not
straight, except
when the height is small compared to $n$ (note however a related result
by \cite{Garet2}, cf. also Remark \ref{rem:Garet2}). This paper aims at filling this gap, although we can do so only in
dimension 2. For instance, suppose that $2h(n)/(nl(A))$
goes to $\tan(\alpha)$ when $n$ goes to infinity, with
$\alpha\in [0,\frac{\pi}{2}]$ and $l(A)$ denoting the length of the
line segment $A$. Our main results imply, under some conditions on $F$ and $A$, that:
\begin{equation}
\label{chapitre6eq:introlgn}
\frac{\phi_n}{nl(A)}\xrightarrow[n\rightarrow \infty]{}\inf_{\widetilde{\theta} \in [\theta-\alpha,\theta+\alpha]} \frac{\nu_{\widetilde{\theta}}}{\cos
  (\widetilde{\theta} - \theta)} \qquad \textrm{a.s. and in }L^1 \;,
\end{equation}
where we re-encoded the function $\nu$ as follows:
$\nu_{\widetilde{\theta}}:=\nu(\vb )$ when
$\vb$ makes an angle $\wt$ with $(1,0)$. Notice that there is no reason for the limit in
(\ref{chapitre6eq:introlgn}) to be identical to $\nu_\theta$. Thus, something different happens when the boxes are not
straight. Notice that this fact can already be observed when $F$ is
concentrated on one point. For instance, if $t(e)=1$
deterministically  and $2h(n)/(nl(A))$
goes to $\tan(\alpha)$ when $n$ goes to infinity, with
$\alpha>\frac{\pi}{4}$, then one may easily compute that
$\nu_\theta=|\cos \theta|+|\sin \theta|$, whereas the limit of
$\phi_n/(nl(A))$ is $\min\{1/|\cos
  \theta |,1/|\sin \theta |\}$. Moreover the moment conditions on $F$ that
  we need to prove (\ref{chapitre6eq:introlgn}) are very weak.

The paper is organized as follows. In section~\ref{chapitre6sec:notations}, we
give the precise definitions and state the main result of the
paper. Section~\ref{chapitre6subsec:concentration} is devoted to a deviation result for $\phi_n$. In
section~\ref{sec:CVesp}, we prove the convergence of the rescaled expectation of
$\phi_n$. Finally, we complete the proof of the law of large numbers for
$\phi_n$ in section~\ref{secLLN}.


\section{Notations,  background and main results}
\label{chapitre6sec:notations}
The most important notations are gathered in sections~\ref{chapitre6subsec:maxflow} to \ref{chapitre6subsec:duality}, the
relevant background is described in section~\ref{chapitre6subsec:mainresults} while our main results are stated in
section~\ref{sec:main}.

\subsection{Maximal flow on a graph}
\label{chapitre6subsec:maxflow}
First, let us define the notion of a flow on a finite
unoriented graph $G=(V,\E)$ with set of vertices $V$ and set of edges
$\E$.  Let
$t=(t(e))_{e\in \E}$ be a collection of non-negative real numbers,
which are called \emph{capacities}. It means that $t(e)$ is the maximal amount of
    fluid that can go through the edge $e$ per unit of time. To each edge $e$, one may associate two oriented edges, and we
shall denote by $\smash{\overrightarrow{\E}}$ the set of all these oriented edges. Let $A$ and $Z$ be
two finite, disjoint, non-empty sets of vertices of $G$: $A$ denotes the source of the
network, and $Z$ the sink. A function $\theta$
on $\smash{\overrightarrow{\E}}$ is called \emph{a flow from $A$ to $Z$ with strength
$\|\theta\|$ and capacities $t$} if it is
antisymmetric, i.e.
$\theta_{\overrightarrow{xy}}=-\theta_{\overrightarrow{yx}}$, if it satisfies the
node law at each vertex $x$ of $V\smallsetminus (A\cup Z)$:
$$\sum_{y\sim x}\theta_{\overrightarrow{xy}}=0\;,$$
where $y\sim x$ means that $y$ and $x$ are neighbours on $G$,
if it satisfies the capacity constraints:
$$\forall e\in \E,\;|\theta(e)|\leq t(e)\;,$$
and if the ``flow in'' at $A$ and the ``flow out'' at $Z$ equal $\|\theta\|$:
$$\|\theta\|=\sum_{a\in A}\sum_{\substack{y\sim a\\ y\not \in
    A}}\theta(\overrightarrow{ay})=\sum_{z\in Z}\sum_{\substack{y\sim
    z\\ y\not \in Z}}\theta(\overrightarrow{yz})\;.$$
The \emph{maximal flow from $A$ to $Z$}, denoted by $\phi_t(G,A,Z)$, is defined as
the maximum strength of all flows from $A$ to $Z$ with capacities
$t$. We shall in general omit the subscript $t$ when it is understood
from the context. The \emph{max-flow min-cut theorem} (see \cite{Bollobas} for instance)
asserts that the maximal flow from $A$ to $Z$ equals the minimal
capacity of a cut between $A$ and $Z$. Precisely, let us say that
$E\subset\E$ is a cut between $A$ and $Z$ in $G$ if every path from
$A$ to $Z$ borrows at least one edge of $E$. Define $V(E)=\sum_{e\in
  E}t(e)$ to be the capacity of a cut $E$. Then,
\begin{equation}
\label{chapitre6eq:maxflowmincut}
\phi_t(G,A,Z)=\min\{V(E)\mbox{ s.t. }E\mbox{ is a cut between
}A\mbox{ and }Z \mbox{ in } G\}\;.
\end{equation}
By convention, if $A$ or $Z$ is empty, we shall define $\phi_t(G,A,Z)$
to be zero.


\subsection{On the square lattice}
\label{chapitre6subsec:squarelattice}
We shall always consider $G$ as a piece of $\ZZ^2$. More precisely, we consider the graph $\LL=(\mathbb{Z}^{2},
\mathbb E ^{2})$ having for vertices $\mathbb Z ^{2}$ and for edges
$\mathbb E ^{2}$, the set of pairs of nearest neighbours for the standard
$L^{1}$ norm. The notation $\langle x,y\rangle$ corresponds to the edge with endpoints
$x$ and $y$. To each edge $e$ in $\mathbb{E}^{2}$ we associate a random
variable $t(e)$ with values in $\mathbb{R}^{+}$. \emph{We suppose that the family
$(t(e), e \in \mathbb{E}^{2})$ is independent and identically distributed,
with a common distribution function $F$}. More
formally, we take the product
measure $\mathbb {P}=F^{\otimes \Omega}$ on $\Omega= \prod_{e\in \mathbb{E}^{2}} [0, \infty[$,
    and we write its expectation $\mathbb{E}$.  If $G$ is
    a subgraph of $\LL$, and $A$ and $Z$ are two subsets of vertices of
    $G$, \emph{we shall denote by $\phi(G,A,Z)$ the maximal flow in $G$ from $A$ to
    $Z$}, where $G$ is equipped with capacities $t$. When $B$ is a
    subset of $\RR^2$, and $A$ and $Z$ are subsets of $\ZZ^2\cap B$,
    we shall denote by $\phi(B,A,Z)$ again the maximal flow
    $\phi(G,A,Z)$ where $G$ is the induced subgraph of $\ZZ^2$ with set of vertices
    $\ZZ^2\cap B$.


We denote by
$\overrightarrow{e}_1$ (resp. $\overrightarrow{e}_2$) the vector $(1,0)\in\RR^2$ (resp. $(0,1)$). Let $A$ be a non-empty
line segment in $\RR^2$. We shall denote by $l(A)$ its (euclidean)
length. All line segments will be
supposed to be closed in $\mathbb{R}^2$. We denote by
$\va$ the vector of unit euclidean norm orthogonal to $\hyp (A)$,
the hyperplane spanned by $A$, and such that there is
$\theta\in[0,\pi[$ such that $\va=(\cos \theta,\sin\theta)$. Define $\vc=( \sin\theta,-\cos\theta)$ and
denote by $a$ and $b$ the end-points of $A$ such that
$(b-a).\vc >0$. For $h$ a positive real number, \emph{we denote by $\cyl(A,h)$ the cylinder of basis $A$ and height $2h$},
i.e., the set 
$$ \cyl (A,h) \,=\, \{x+t \va \,|\, x\in A \,,\,  t\in
[-h,h]    \}\,.$$
We define also \emph{the $r$-neighbourhood $\mathcal{V} (H,r)$ of a subset $H$ of $\mathbb{R}^d$} as
$$\mathcal{V}(H,r) \,=\, \{ x \in \mathbb{R}^d \,|\, d(x,H)<r\}\,,$$
where the distance is the euclidean one ($d(x,H) = \inf \{\|x-y\|_2
\,|\, y\in H \}$).

\begin{figure}[!ht]
\centering
\begin{picture}(0,0)%
\includegraphics{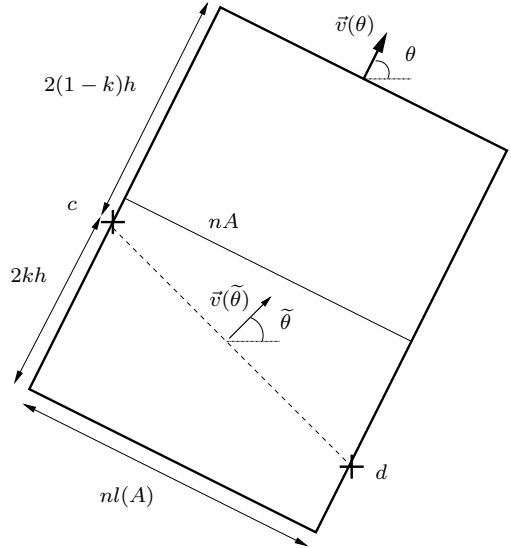}%
\end{picture}%
\setlength{\unitlength}{1973sp}%
\begingroup\makeatletter\ifx\SetFigFont\undefined%
\gdef\SetFigFont#1#2#3#4#5{%
  \reset@font\fontsize{#1}{#2pt}%
  \fontfamily{#3}\fontseries{#4}\fontshape{#5}%
  \selectfont}%
\fi\endgroup%
\begin{picture}(6240,6879)(1294,-8053)
\put(2101,-3811){\makebox(0,0)[rb]{\smash{{\SetFigFont{8}{9.6}{\rmdefault}{\mddefault}{\updefault}{\color[rgb]{0,0,0}$c$}%
}}}}
\put(5851,-7186){\makebox(0,0)[lb]{\smash{{\SetFigFont{8}{9.6}{\rmdefault}{\mddefault}{\updefault}{\color[rgb]{0,0,0}$d$}%
}}}}
\put(6226,-1936){\makebox(0,0)[lb]{\smash{{\SetFigFont{8}{9.6}{\rmdefault}{\mddefault}{\updefault}{\color[rgb]{0,0,0}$\theta$}%
}}}}
\put(4651,-5311){\makebox(0,0)[lb]{\smash{{\SetFigFont{8}{9.6}{\rmdefault}{\mddefault}{\updefault}{\color[rgb]{0,0,0}$\wt$}%
}}}}
\put(3751,-4036){\makebox(0,0)[lb]{\smash{{\SetFigFont{8}{9.6}{\rmdefault}{\mddefault}{\updefault}{\color[rgb]{0,0,0}$nA$}%
}}}}
\put(4276,-5011){\makebox(0,0)[rb]{\smash{{\SetFigFont{8}{9.6}{\rmdefault}{\mddefault}{\updefault}{\color[rgb]{0,0,0}$\vec{v}(\wt)$}%
}}}}
\put(5851,-1561){\makebox(0,0)[rb]{\smash{{\SetFigFont{8}{9.6}{\rmdefault}{\mddefault}{\updefault}{\color[rgb]{0,0,0}$\vec{v}(\theta)$}%
}}}}
\put(2851,-2311){\makebox(0,0)[rb]{\smash{{\SetFigFont{8}{9.6}{\rmdefault}{\mddefault}{\updefault}{\color[rgb]{0,0,0}$2(1-k)h$}%
}}}}
\put(1726,-4711){\makebox(0,0)[rb]{\smash{{\SetFigFont{8}{9.6}{\rmdefault}{\mddefault}{\updefault}{\color[rgb]{0,0,0}$2kh$}%
}}}}
\put(3076,-7486){\makebox(0,0)[rb]{\smash{{\SetFigFont{8}{9.6}{\rmdefault}{\mddefault}{\updefault}{\color[rgb]{0,0,0}$nl(A)$}%
}}}}
\end{picture}%
\caption{An admissible boundary condition $(k,\wt)$.}
\label{chapitre6fig:notations}
\end{figure}

Now, $D(A,h)$ denotes the set of \emph{admissible boundary conditions} on
$\cyl(A,h)$ (see Figure \ref{chapitre6fig:notations}): 
$$D(A,h)=\left\{(k,\tilde \theta) \,|\, k\in [0,1]\mbox{ and
  }\tilde\theta\in
  \left[\theta-\arctan\left(\frac{2hk}{l(A)}\right),\theta+\arctan\left(\frac{2h(1-k)}{l(A)}\right)\right]\right\}\;.$$

The meaning of an element $\kappa=(k,\tilde \theta)$ of $D(A,h)$ is the
following. We define
$$\vb \,=\, (\cos \widetilde{\theta}, \sin
\widetilde{\theta}) \qquad \textrm{and} \qquad \vd
\,=\, (\sin \widetilde{\theta}, -\cos \widetilde{\theta})\,. $$
In $\cyl(nA,h(n))$, we may define two points $c$ and $d$
such that $c$ is ``at height $2kh$ on the left side of $\cyl(A,h)$'',
and $d$ is ``on the right side of $\cyl(A,h)$''  by 
$$c=a+(2k-1)h\va \,,\quad
(d-c)\mbox{ is orthogonal to }\vb \quad\mbox{ and
} \quad d\mbox{ satisfies }\vec{cd}\cdot \vd >0 \;.$$
Then we see that $D(A,h)$ is exactly the set of parameters so that $c$ and $d$
remain ``on the sides of $\cyl(A,h)$''. 

We define also $\D(A,h)$, the
set of angles $\tilde\theta$ such that there is an admissible boundary
condition with angle $\tilde\theta$:
$$\D(A,h)=\left[\theta-\arctan\left(\frac{2h}{l(A)}\right),\theta+\arctan\left(\frac{2h}{l(A)}\right)\right]\;.$$
It will be useful to define
 the \emph{left side (resp. right side) of $\cyl(A,h)$}: let $\lleft(A)$
 (resp. $\rright(A)$) be the set of vertices in $\cyl(A,h)\cap\ZZ^2$
such that there exists $y\notin \cyl(A,h)$, $\langle x,y\rangle \in
\mathbb{E}^d$ and $[ x,y [$, the segment that includes $x$ and excludes
$y$, intersects $a+[-h,h].\va $ (resp. $b+[-h,h].\va $).

Now, the set $\cyl(A,h) \smallsetminus (c+\RR(d-c))$ has two connected
components, which we denote by $\mathcal{C}_1(A,h,k,\tilde\theta)$ and
$\mathcal{C}_2(A,h,k,\tilde\theta)$. For $i=1,2$, let $A_i^{h,k,\tilde\theta}$ be
the set of the points in $\mathcal{C}_i(A,h,k,\tilde\theta) \cap \mathbb{Z}^2$ which have
a nearest neighbour in $\mathbb{Z}^2 \smallsetminus \cyl(A,h)$:
$$A_i^{h,k,\tilde\theta}\,=\,\{x\in \mathcal{C}_i(A,h,k,\tilde\theta) \cap
\mathbb{Z}^2 \,|\, \exists y \in \mathbb{Z}^2 \smallsetminus \cyl(A,h) \,,\,
\|x-y\|_{1} =1 \}\,.$$

We define \emph{the flow in $\cyl(A,h)$ constrained by the boundary
condition $\kappa=(k,\tilde\theta)$} as:
$$\phi^\kappa(A,h):=\phi(\cyl(A,h),A_1^{h,k,\tilde\theta},A_2^{h,k,\tilde\theta})\;.$$
A special role is played by the condition  $\kappa = (1/2, \theta)$,
and we shall denote: 
$$\tau(A,h) = \tau (\cyl(A,h), \va )=\phi^{(1/2, \theta)}(A,h)\;.$$
Let $T(A,h)$ (respectively $B(A,h)$) be the top
(respectively the bottom) of $\cyl(A,h)$, i.e.,
$$ T(A,h) \,=\, \{ x\in \cyl(A,h) \,|\, \exists y\notin \cyl(A,h)\,,\,\,
\langle x,y\rangle \in \mathbb{E}^d \mbox{ and }\langle x,y\rangle
\mbox{ intersects } A+h\va  \}  $$
and
$$  B(A,h) \,=\, \{ x\in \cyl(A,h) \,|\, \exists y\notin \cyl(A,h)\,,\,\,
\langle x,y\rangle \in \mathbb{E}^d \mbox{ and }\langle x,y\rangle
\mbox{ intersects } A-h\va \} \,.$$
We shall denote the flow in $\cyl(A,h)$ from the top to the bottom as:
$$\phi(A,h)= \phi (\cyl(A,h), \va)=\phi(\cyl(A,h),T(A,h),B(A,h))\;.$$


\subsection{Duality}
\label{chapitre6subsec:duality}

The main reason why dimension 2 is easier to deal with than dimension
$d\geq 3$ is duality. Planar duality implies that there are only
$O(h^2)$ admissible boundary conditions on $\cyl(A,h)$. Let us go a
bit into the details.

The dual lattice $\LL^*$ of $\LL$ is constructed as follows:
place a vertex in the centre of each face of $\LL$ and join two
vertices in $\LL^*$ if and only if the corresponding faces of $\LL$ share
an edge. To each edge $e^*$ of $\LL^*$, we assign the time coordinate
$t(e)$, where $e$ is the unique edge of $\mathbb{E}^2$ crossed by
$e^*$. 
Now, let $A$ be a line segment in $\RR^2$. Let $G_A$ be the induced
subgraph of $\LL$ with set of vertices $\cyl(A,h)\cap \ZZ^2$. Let
$G_A^*$ be the planar dual of $G_A$ in the following sense: $G_A^*$ has set
of edges $\{e^*\mbox{ s.t. }e\in G_A\}$, and set of vertices those
vertices which belong to  this set of edges. Now, we define $\lleft^*(A)$
(resp. $\rright^*(A)$) as the set of vertices $v$
  of $G_A^*$ which have at least one neighbour in $\LL^*$ which is not
  in $G_A$ and such that there exists an edge $e^*$ in $G_A^*$ with
  $v\in e^*$ and $e^*\cap \lleft(A)\not =\emptyset$ (resp. $e^*\cap \rright(A)\not =\emptyset$).

It is well known that the (planar) dual of a cut between the top and the
bottom of $\cyl(A,h)$ is a self-avoiding path from ``left'' to
``right''. Furthermore, if the cut is minimal for the inclusion,
the dual self-avoiding path has only one vertex on the left boundary
of the dual of $A\cap\ZZ^2$ and one vertex on the right boundary. The
following lemma is a formulation in our setting of those classical duality
results (see for instance \cite{GrimmettKesten84} p.358 and \cite{Bollobas},
p.47).
\begin{lem}
Let $A$ be a line segment $\RR^2$ and $h$ be a positive real
number. If $E$ is a set of edges, let 
$$E^*=\{e^* \,|\,e\in E\}\;.$$
If $E$ is a cut between $B(A,h)$ and $T(A,h)$, minimal for the
inclusion, then $E^*$ is a self-avoiding path from $\lleft^*(A)$ to
$\rright^*(A)$ such that exactly one point of $E^*$ belongs to
$\lleft^*(A)$, exactly one point of $E^*$ belongs to
$\rright^*(A)$, and these two points are the end-points of the path.
\end{lem}
An immediate consequence of this planar duality is the following.
\begin{lem}
\label{chapitre6lem:duality}
Let $A$ be any line segment in $\RR^2$ and $h$ a positive real
number. Then,
\begin{equation*}
\phi(A,h)= \min_{\kappa\in D(A,h)}\phi^\kappa(A,h) \;.
\end{equation*}
\end{lem}
Notice that the condition $\kappa$ belongs to the non-countable set
$D(A,h)$, but the graph is discrete so $\phi^\kappa(A,h)$ takes only a
finite number of values when $\kappa\in D(A,h)$. Precisely, there is a
finite subset $\tilde D(A,h)$ of $D(A,h)$, such that:
\begin{equation}
\label{chapitre6eq:Oh2}
\card(\tilde D(A,h))\leq C_4 h^2\;,
\end{equation}
for some universal constant $C_4$, and:
\begin{equation*}
\phi(A,h)= \min_{\kappa\in \tilde D(A,h)}\phi^\kappa(A,h) \;.
\end{equation*}


\subsection{Background}
\label{chapitre6subsec:mainresults}

First, let us recall some facts concerning the behaviour of
$\tau(nA,h(n))$ when $n$ and $h(n)$ go to infinity. Using a subadditive argument and deviation inequalities, Rossignol
and Th\'eret have proved in \cite{RossignolTheret08b} that $\tau(nA,
h(n))$ satisfies a law of large numbers:
\begin{thm}
\label{chapitre6thm:LGNtau}
We suppose that
$$ \int_{[0,\infty[} x \, dF (x) \,<\, \infty\, .  $$
For every unit vector $\va =(\cos\theta,\sin\theta)$, there exists a
constant $\nu_\theta$ depending on $F$, $d$ and $\theta$, such that for
every non-empty line-segment $A$ orthogonal to $\va$ and of euclidean length $l(A)$, for
every height function $h: \NN \rightarrow \RR^+$ satisfying
$\lim_{n\rightarrow \infty} h(n) = +\infty$, we have
$$ \lim_{n\rightarrow \infty} \frac{\tau(nA, h(n))}{n l(A)}
\,=\, \nu_\theta \qquad \textrm{in } L^1 \,. $$
Moreover, if the origin of the graph belongs to $A$, or if
$$ \int_{[0,\infty[} x^{2} \, dF (x) \,<\, \infty \,, $$
then
$$ \lim_{n\rightarrow \infty} \frac{\tau(nA, h(n))}{nl(A)}
\,=\, \nu_\theta \qquad \textrm{a.s.} $$
Under the added assumption that $\lim_{n\rightarrow \infty} h(n) /n =0$,
the variable $\phi(nA,h(n))$ satisfies the same law of large numbers as
$\tau(nA,h(n))$, under the same conditions.
\end{thm}
This law of large numbers holds in fact for every dimension $d\geq 2$. Let
us remark that (in dimension two) $\nu_\theta$ is equal to $\mu(\vc ) = \mu(\va)$, where
$\mu(.)$ is the time-constant function of first passage percolation as
defined in \cite{Kesten:StFlour}, (3.10) p.~158. This equality follows
from the duality considerations of section \ref{chapitre6subsec:duality} and
standard first passage percolation techniques (see also
Theorem 5.1 in \cite{GrimmettKesten84}) that relate cylinder passage times to
unrestricted passage times (as in \cite{HammersleyWelsh},
Theorem~4.3.7 for instance). Boivin has also proved a very similar law of large
numbers (see Theorem 6.1 in \cite{Boivin}). Notice that for the
definition of $\mu(.)$, Kesten requires only the existence of the first
moment of the law $F$ in the proof from \cite{Kesten:StFlour}, and it
can also be defined under the weaker condition $\int_0^{\infty}(1
-F(x))^4\;dx<\infty$. 

One consequence of this equality between $\nu$ and $\mu$ is that
$\theta\mapsto \nu_\theta$ is either constant equal to zero, or always
non-zero. In fact the following property holds (cf. \cite{Kesten:StFlour},
Theorem 6.1 and Remark 6.2 p.~218):
\begin{prop}
\label{propnu}
We suppose that $\int_{[0,+\infty[} x \, dF(x) <\infty$. Then $\nu_\theta$
is well defined for all $\theta$, and we have
$$ \nu_\theta \,>\,0 \iff F(0) \,<\, 1/2\,.  $$
\end{prop}

There exists a law of large numbers for the variable $\phi(nA, h(n))$ when
the rectangle we consider is straight, i.e., $\theta = 0$. It has been
proved in \cite{GrimmettKesten84}, Corollary 4.2, that:
\begin{thm}
\label{thmllnphidroit}
Suppose that $A=[0,1]\times\{0\}$, $\int_{[0,+\infty[} x \, dF(x) <\infty$,
$$h(n)\xrightarrow[n\rightarrow\infty]{}\infty \quad\mbox{ and
}\quad\frac{\log h(n)}{n}\xrightarrow[n\rightarrow\infty]{}0\;.$$
Then,
$$\frac{\phi(nA, h(n))}{n}\xrightarrow[n\rightarrow\infty]{a.s}\nu_0\;.$$
\end{thm}
\begin{rem}
Notice that in \cite{GrimmettKesten84}, the condition on $F$ is in
fact weakened to $\int_0^{\infty}(1 -F(x))^4\;dx<\infty$, obtaining
the convergence to $\mu((0,1))$. However, our definition of $\nu_\theta$ requires a moment of order 1. 
\end{rem}

Finally, let us remark that \cite{Garet2} proved a law of large
numbers for the maximal flow between a compact convex set
$\Lambda \subset\RR^2$ and infinity. This is somewhat related to our main
result, Theorem \ref{chapitre6thm:lgn}, see Remark \ref{rem:Garet2}.
Before stating Garet's result, we need some notations. For every convex
bounded set $\Lambda \subset \RR^2$, we denote by $\partial ^* \Lambda$ the
set of all the points $x$ of the boundary $\partial \Lambda$ of $\Lambda$
where $\Lambda$ admits a unique outer normal, which is denoted by
$\vec{v}_{\Lambda} (x)$. We denote the coordinates of $\vec{v}_{\Lambda}
(x)$ by $(\cos(\theta(\Lambda,x)), \sin(\theta(\Lambda, x)))$ for every $x$ in $\partial ^* \Lambda$. We denote by $\sigma (\Lambda)$ the maximal flow
from $\Lambda$ to infinity. Let $\H^1$ be the one-dimensional Hausdorff
measure. Theorem 2.1 in \cite{Garet2} is the following:
\begin{thm}
\label{thmgaret}
We suppose that $F(0)<1/2$ and that
$$ \exists \gamma>0 \qquad \int_{[0,+\infty[} e^{\gamma t(e)} \,dF(x)
\,<\,\infty \,.  $$
Then, for each bounded convex set $\Lambda \subset \RR^2$ with the origin
of the graph $0$ in its interior, we have
\begin{equation}
\label{eqgaret}
\lim_{n\rightarrow \infty} \frac{\sigma(n\Lambda)}{n} \,=\, \int_{\partial
  ^* \Lambda} \nu_{\theta(\Lambda, x)} d\H^1(x) \,=\, \mathcal{I} (\Lambda)
\,>\,0\,.
\end{equation}
\end{thm}


\subsection{Main result}
\label{sec:main}
We recall that for all $n\in \NN$, we have defined
$$\D(nA,h(n))=\left[\theta-\arctan\left(\frac{2h(n)}{nl(A)}\right),\theta+\arctan\left(\frac{2h(n)}{nl(A)}\right)\right]\;.$$
We may now state our main result.
\begin{thm}
\label{chapitre6thm:lgn}
Let $A$ be a non-empty line-segment in $\RR^2$, with euclidean length
$l(A)$. Let $\theta\in[0,\pi[$ be such that $(\cos\theta,\sin\theta)$
is orthogonal to $A$ and $(h(n))_{n\geq 0}$
be a sequence of positive real numbers such that:
\begin{equation}
\label{chapitre6eq:conditionshn}
\left\lbrace \begin{array}{l}h(n)\xrightarrow[n\rightarrow
    \infty]{}+\infty\;,\\
\frac{\log h(n)}{n}\xrightarrow[n\rightarrow
    \infty]{}0\;.\end{array}\right.
\end{equation}
Define:
$$\overline{\D}=\limsup_{n\rightarrow
  \infty}\D(nA,h(n))=\bigcap_{N\geq 1}\bigcup_{n\geq
  N}\D(nA,h(n))\;,$$
and
$$\underline{\D}=\liminf_{n\rightarrow \infty}\D(nA,h(n))=\bigcup_{N\geq 1}\bigcap_{n\geq N}\D(nA,h(n))\;.$$
Suppose that $F$ has a finite moment of order 1:
\begin{equation}
\label{chapitre6eq:conditionsFLGN}
\int_0^\infty x \, dF(x)<\infty\;.
\end{equation}
Then,
\begin{equation}
\label{eq:liminfmoyenne}
\liminf_{n\rightarrow \infty}\frac{\EE[\phi(nA,h(n))]}{nl(A)}=\inf \left\{
  \frac{\nu_{\widetilde{\theta}}}{\cos  (\widetilde{\theta} - \theta)}
  \,|\, \widetilde{\theta} \in \overline{\D} \right\}
\end{equation}
and
\begin{equation}
\label{eq:limsupmoyenne}\limsup_{n\rightarrow\infty}\frac{\EE[\phi(nA,h(n))]}{nl(A)}= \inf\left\{
  \frac{\nu_{\widetilde{\theta}}}{\cos (\widetilde{\theta} - \theta)} \,|\,
  \widetilde{\theta} \in \underline{\D}  \right\} \,.
\end{equation}
Moreover, if
\begin{equation}
\label{chapitre6eq:conditionsFLGNbis}
\int_0^\infty x^2 \, dF(x)<\infty\,,
\end{equation}
or if:
\begin{equation}
\label{conditionO}
0\mbox{ is the middle of }A\;,
\end{equation}
then
$$\liminf_{n\rightarrow \infty}\frac{\phi(nA,h(n))}{nl(A)}=\inf \left\{
  \frac{\nu_{\widetilde{\theta}}}{\cos  (\widetilde{\theta} - \theta)}
  \,|\, \widetilde{\theta} \in \overline{\D} \right\} \qquad a.s.$$
and
$$\limsup_{n\rightarrow\infty}\frac{\phi(nA,h(n))}{nl(A)}= \inf\left\{
  \frac{\nu_{\widetilde{\theta}}}{\cos (\widetilde{\theta} - \theta)} \,|\,
  \widetilde{\theta} \in \underline{\D}  \right\} \qquad a.s.$$
\end{thm}

\begin{rem} It is likely that condition
  (\ref{chapitre6eq:conditionsFLGN}) can be weakened to $\int_0^\infty
  (1-F(x))^4\;dx$, as in Theorem \ref{thmllnphidroit}. This would
  require to define $\nu$ a bit differently.
\end{rem}

\begin{cor}
\label{chapitre6corollaire}
We suppose that conditions (\ref{chapitre6eq:conditionshn}) on $h$ are
satisfied. We suppose also that there is some
$\alpha\in \left[0,\frac{\pi}{2}\right]$ such that:
$$\frac{2h(n)}{nl(A)}\xrightarrow[n\rightarrow
    \infty]{}\tan \alpha \;.$$
Then, if condition (\ref{chapitre6eq:conditionsFLGN}) on $F$ is
satisfied, we have
$$\lim_{n\rightarrow \infty}\frac{\phi(nA,h(n))}{nl(A)}=\inf\left\{
  \frac{\nu_{\widetilde{\theta}}}{\cos (\widetilde{\theta} - \theta)} \,|\,
  \widetilde{\theta} \in [\theta-\alpha,\theta+\alpha] \right\} \qquad
\textrm{in } L^1\,.$$
Moreover, if condition (\ref{chapitre6eq:conditionsFLGNbis}) or
(\ref{conditionO}) are satisfied, then
$$\lim_{n\rightarrow \infty}\frac{\phi(nA,h(n))}{nl(A)}=\inf\left\{
  \frac{\nu_{\widetilde{\theta}}}{\cos (\widetilde{\theta} - \theta)} \,|\,
  \widetilde{\theta} \in [\theta-\alpha,\theta+\alpha] \right\} \qquad
\textrm{a.s.}$$
\end{cor}
It has already been remarked in \cite{Theret:small} (see the
discussion after Theorem 2) that the condition on $h$ is the good one
to have positive speed when one allows edge capacities to be
null with positive probability. 

\begin{rem}
Notice that Theorem \ref{chapitre6thm:lgn} is consistent with Theorem
\ref{thmllnphidroit}, the existing law of large numbers for $\phi(nA, h(n)$
in the straight case. Indeed, it is
known that $\nu$ satisfies the weak triangle inequality (see section 4.4 in \cite{RossignolTheret08b}), and for
symmetry reasons, it implies that when $\theta \in
\{0,\pi /2 \}$, the function $\wt\mapsto \nu_{\widetilde{\theta}} / \cos
  (\widetilde{\theta} - \theta)$ is minimum for $\wt=\theta$ and thus, Theorem \ref{chapitre6thm:lgn} implies that
$\phi(nA,h(n))/(nl(A))$ converges to $\nu_0$, the limit of
$\tau(nA,h(n))/(nl(A))$, when $\cyl (nA, h(n))$ is a straight
cylinder. In fact, the same phenomenon occurs for any $\theta$ such
that there is a symmetry axis of direction $\theta$ for the
lattice $\ZZ^2$. These directions in $[0,\pi[$ are
of course $\{0, \pi /4, \pi/ 2, 3 \pi /4\}$. Also, Corollary \ref{chapitre6corollaire} is consistent with the fact that for
general boxes, when $h(n)$ is small with respect to $n$,
$\phi(nA,h(n))/(nl(A))$ and $\tau(nA,h(n))/(nl(A))$ have the same
limit. 
\end{rem}

\begin{rem}
\label{rem:Garet2}
Theorems \ref{chapitre6thm:lgn} and \ref{thmgaret} are
related. First, they are stated in dimension two only, because both proofs
use the duality of the planar graph to define the path which is the dual of a cutset, and then
the fact that such paths can be glued together if they have a common
endpoint. These properties hold only in dimension two: the dual of
an edge in dimension greater than three is a unit surface, and it is much
more difficult to study the boundary of a surface. This is the reason why
these theorems are not yet generalized in higher dimensions (see also Remark
\ref{remd3}). Moreover, the expressions of the limits $\mathcal{I}(\Lambda)$
and $\eta_{\theta,h}$ appearing in these theorems are very similar. On one
hand, the
constant $\eta_{\theta,h}$ is the infimum of the integral of $\nu$ along
the segments that cut the top from the bottom of $\cyl(A,h(n)/n)$ for large
$n$. Since $\nu$ satisfies the weak triangle inequality, $\eta_{\theta, h}$
is also equal to infimum of the integral of $\nu$ along the polyhedral
curves that have the same property of cutting. On the other hand,
Garet only has  to consider the case of a polyhedral convex set $\Lambda$ during his
proof, and he proves the important following property: if $\Lambda \subset
\Lambda '$, where $\Lambda$ and $\Lambda'$ are polyhedral and $\Lambda$ is
convex, then $\mathcal{I}(\Lambda) \leq \mathcal{I}(\Lambda ')$. Thus, for
a polyhedral convex set $\Lambda$, $\mathcal{I}(\Lambda)$ is the infimum of
the integral of $\nu$ along the polyhedral curves that cut $\Lambda$ from
infinity.
\end{rem}



\subsection{Sketch of the proof}
\label{chapitre6subsec:sketchLGN}

We suppose that $A$ is a non-empty line segment in
$\RR^2$. To shorten the
notations, we shall write $D_n=D(nA,h(n))$, the set of all admissible conditions for $(nA,h(n))$:
$$D_n=\left\{(k,\tilde \theta) \,|\,k\in [0,1]\mbox{ and
   }\tilde\theta\in
  \left[\theta-\arctan\left(\frac{2h(n)k}{nl(A)}\right),\theta+\arctan\left(\frac{2h(n)(1-k)}{nl(A)}\right)\right]\right\}\;,$$
and
$$\D_n=\left[\theta-\arctan\left(\frac{2h(n)}{nl(A)}\right),\theta+\arctan\left(\frac{2h(n)}{nl(A)}\right)\right]\;.$$
Also, we shall use:
$$\phi_n=\phi(nA,h(n)),\quad
\phi_n^\kappa=\phi^\kappa(nA,h(n))\quad\mbox{ and }\quad\tau_n=\tau(nA,h(n))\;.$$

First, notice that $0\leq\phi_n\leq\tau_n$. If $F(0)\geq
\frac{1}{2}$, then $\tau_n/n$ converges to zero, and so does $\phi_n$,
so Theorem \ref{chapitre6thm:lgn} is trivially true. {\it We shall therefore make the
  following hypothesis in the rest of the article:}
\begin{equation}
\label{eq:F(0)}
F(0)<\frac{1}{2}\;.
\end{equation}

Now, let us draw a sketch of the proof of Theorem \ref{chapitre6thm:lgn}. Recall that from Lemma
\ref{chapitre6lem:duality},
$$\phi_n=\min_{\kappa\in D_n}\phi_n^\kappa\;.$$
First, we shall study the asymptotics of $\EE(\phi_n)$ (section 4):

{\bf Step 1.} By a subadditive argument (see
Figure \ref{chapitre6emboitement1}), we show in section
\ref{chapitre6subsec:upperbound} that
$$\limsup_{n\rightarrow \infty} \frac{\EE[\phi_n]}{nl(A)} \,\leq\, \inf_{\wt \in
\underline{\D}} \frac{\nu_{\wt}}{\cos (\wt - \theta)}  $$
and
$$\liminf_{n\rightarrow \infty} \frac{\EE[\phi_n]}{nl(A)} \,\leq\, \inf_{\wt \in
\overline{\D}} \frac{\nu_{\wt}}{\cos (\wt - \theta)} \,. $$

{\bf Step 2.} On the other hand, by a similar subadditive argument (see Figure
\ref{chapitre6emboitement2}), we show in section
\ref{chapitre6subsec:lowerbound} that
$$\liminf_{n\rightarrow \infty}\inf_{\kappa \in D_n}
\frac{\EE[\phi_n^{\kappa}]}{nl(A)} \,\geq\,  \inf_{\wt \in
\overline{\D}} \frac{\nu_{\wt}}{\cos (\wt - \theta)}  $$
and
$$\limsup_{n\rightarrow \infty}\inf_{\kappa \in D_n}
\frac{\EE[\phi_n^{\kappa}]}{nl(A)} \,\geq\,  \inf_{\wt \in
\underline{\D}} \frac{\nu_{\wt}}{\cos (\wt - \theta)} \,. $$

{\bf Step 3.} Using deviation results for the variables
$\phi_n^\kappa$ (section \ref{chapitre6subsec:concentration}),  we prove in section \ref{chapitre6subsec:endLGN} that $\EE[\phi_n]$ is equivalent to
$\inf_{\kappa \in D_n} \EE[\phi_n^{\kappa}]$, and this ends the study of the
asymptotic behaviour of $\EE[\phi_n]$.

Next, we relate $\phi_n$ and $\EE(\phi_n)$ to show the almost sure
asymptotics (section \ref{secLLN}):

{\bf Step 4.} A deviation result for $\phi_n$ obtained in section
\ref{chapitre6subsec:concentration} shows that almost surely, asymptotically, $\phi_n/n$
is at least as large as $\EE(\phi_n)/n$. 

{\bf Step 5.} Finally, we use again the subadditive argument of the
first step of the proof to prove that almost surely, 
$$\limsup_{n\rightarrow \infty} \frac{\phi_n}{nl(A)} \,\leq\, \inf_{\wt \in
\underline{\D}} \frac{\nu_{\wt}}{\cos (\wt - \theta)}  $$
and
$$\liminf_{n\rightarrow \infty} \frac{\phi_n}{nl(A)} \,\leq\, \inf_{\wt \in
\overline{\D}} \frac{\nu_{\wt}}{\cos (\wt - \theta)} \,. $$

\section{Deviation properties of the maximal flows}
\label{chapitre6subsec:concentration}

The following proposition, due to Kesten, allows to control the size
of the minimal cut, and is of fundamental importance in the study of
First Passage Percolation. 
\begin{prop}[Proposition 5.8 in \cite{Kesten:StFlour}]
\label{chapitre6prop:5.8Kesten}
Suppose that $F(0)<\frac{1}{2}$. Then, there are constants $\eps$, $C_1$
and $C_2$, depending only on $F$, such that:
$$\PP\left( \begin{array}{c} \exists\mbox{ a self-avoiding path }\gamma\mbox{ in }\LL^*\mbox{,
  starting at }(\frac{1}{2},\frac{1}{2})\mbox{,} \\ \mbox{with }\card(\gamma)\geq m\mbox{ and
}\sum_{e^*\in\gamma}t(e^*)\leq \eps m \end{array} \right) \,\leq\, C_1e^{-C_2m}\;.$$
\end{prop}
Thanks to Proposition \ref{chapitre6prop:5.8Kesten} and general deviation
 inequalities due to \cite{Boucheronetal03}, we obtain the
following deviation result for the maximal flows $\phi_n$ and
$\phi_n^\kappa$. The proof is exactly the same as the proof of
Proposition~4.3 in \cite{RossignolTheret08b}, using Proposition
\ref{chapitre6prop:5.8Kesten} instead of Zhang's result. We reproduce
it here for the sake of completeness.
\begin{prop}
\label{prop:deviation}
Suppose that hypotheses (\ref{chapitre6eq:conditionsFLGN}) and
(\ref{eq:F(0)}) hold. Then,
for any $\eta\in]0,1]$, there are strictly positive constants $C(\eta,F)$,
$K_1(F)$ and $K_2(F)$, such that, for every $n\in\NN^*$, and every
non-degenerate line segment $A$,
\begin{equation}
\label{chapitre6eq:deviationkappa}
 \max_{\kappa\in
  D_n}\PP(\phi_n^\kappa<\EE(\phi_n^\kappa)(1-\eta))\leq
K_1e^{-C(\eta,F)\min_\kappa\EE(\phi_n^\kappa)}\;.
\end{equation}
and:
\begin{equation}
\label{eq:dev_phi_ordre1}
\PP\left( \phi_n \leq \EE(\phi_n)(1-\eta) \right) \leq K_2h(n)^2e^{-C(\eta,F)\EE(\phi_n)}\;.
\end{equation}
\end{prop}
\begin{dem}
Let us fix $A$, $n\in\NN^*$ and $\kappa=(k,\tilde\theta)\in D_n$. First, we prove the
result for $\phi_n^\kappa$. We shall denote by $E_{\phi_n^\kappa}$ a cut whose capacity achieves
the minimum in the dual definition (\ref{chapitre6eq:maxflowmincut}) of
$\phi_n^\kappa$.  Since
$\PP\left( \phi_n^\kappa\leq \EE(\phi_n^\kappa)(1-\eta) \right)$ is a decreasing
function of $\eta$, it is enough to prove the result for all $\eta$
less than or equal to some absolute $\eta_0\in]0,1[$. We use this remark to
exclude the case $\eta =1$ in our study, thus, from now on,
let $\eta$ be a fixed real number in $]0,1[$.

We order the edges in $\cyl(nA,h(n))$ as
$e_1,\ldots,e_{m_n}$. For every hyperrectangle $A$, we denote by
$\mathcal{N}(A,h)$ the minimal number of edges in $A$ that can
disconnect $A_1^h$ from $A_2^h$ in $\cyl(A,h)$. For any real number $r\geq \mathcal{N}(nA,h(n))$, we define:
$$ \psi_n^r\,=\, \min\left\{ \begin{array}{c} V(E)\mbox{ s.t. }\card(E)\leq
r\mbox{ and } E \mbox{ cuts }\\
(nA)_1^{h(n),k,\tilde\theta}\mbox{
  from }(nA)_2^{h(n),k,\tilde\theta}\mbox{ in }\cyl(nA,h(n)) \end{array} \right\}\;. $$
Now, suppose that hypotheses (\ref{chapitre6eq:conditionsFLGN}) and
(\ref{eq:F(0)}) hold, let $\eps$, $C_1$ and $C_2$ be as
in Proposition~\ref{chapitre6prop:5.8Kesten}, and define
$r=(1-\eta)\EE(\phi_n^\kappa)/\eps$.  Suppose
first that $r<\mathcal{N}(nA,h(n))$. Then,
\begin{eqnarray*}
\PP(\phi_n^\kappa\leq (1-\eta)\EE(\phi_n^\kappa))&=&\PP(\phi_n^\kappa\leq
(1-\eta)\EE(\phi_n^\kappa)\mbox{ and } \card (E_{\phi_n^\kappa}) \geq
(1-\eta)\EE(\phi_n^\kappa)/\eps)\;,\\
&\leq&C_1e^{-C_2(1-\eta)\EE(\phi_n^\kappa)/\eps}\;,
\end{eqnarray*} 
from Proposition~\ref{chapitre6prop:5.8Kesten}, and the desired inequality is
obtained. Suppose now that we have $r\geq\mathcal{N}(nA,h(n))$. Then,
\begin{eqnarray}
\nonumber \PP(\phi_n^\kappa\leq (1-\eta)\EE(\phi_n^\kappa))&=&\PP(\phi_n^\kappa\leq
(1-\eta)\EE(\phi_n^\kappa)\mbox{ and }\psi_n^r\not=\phi_n^\kappa)+\PP(\psi_n^r\leq (1-\eta)\EE(\phi_n^\kappa))\;,\\
\label{eq:tautaunr} &\leq&C_1e^{-C_2r}+\PP(\psi_n^r\leq (1-\eta)\EE(\psi_n^r))\;,
\end{eqnarray}
from Proposition~\ref{chapitre6prop:5.8Kesten} and the fact that $
\phi_n^\kappa\leq \psi_n^r$. Now, we truncate our variables $t(e)$. Let $a$ be a positive
real number to be chosen later, and define $\tilde t(e)=t(e)\land
a$. Let:
$$ \tilde \psi_n^r\,=\, \min\left\{ \begin{array}{c} \sum_{e\in
      E}\tilde t(e)\mbox{ s.t. }\card(E)\leq r\mbox{ and }E \mbox{ cuts }\\
 (nA)_1^{h(n),k,\tilde\theta}\mbox{
  from }(nA)_2^{h(n),k,\tilde\theta}\mbox{ in }\cyl(nA,h(n)) \end{array} \right\}\;. $$
Notice that $\tilde \psi_n^r\leq \psi_n^r$. We shall denote by
$E_{\tilde\psi_n^r}$ a cutset whose capacity achieves the minimum in the
definition of $\tilde \psi_n^r$. If there are more than one, we use a
deterministic method to select a unique one with the minimal number of
edges among these. Then,
\begin{eqnarray*}
0\leq \EE(\psi_n^r)-\EE(\tilde\psi_n^r)&\leq &\EE\left[\sum_{e\in
  E_{\tilde\psi_n^r}}t(e)-\sum_{e\in
  E_{\tilde\psi_n^r}}\tilde t(e)\right]\;,\\
&\leq &\EE\left[\sum_{e\in
  E_{\tilde\psi_n^r}}t(e)\II_{t(e)\geq a}\right]\;,\\
&=&\sum_{i=1}^{m_n}\EE(t(e_i)\II_{t(e_i)\geq a}\II_{e_i\in
  E_{\tilde\psi_n^r}})\;,\\
&=&\sum_{i=1}^{m_n}\EE\left\lbrack \EE\left(t(e_i)\II_{t(e_i)\geq a}\II_{e_i\in
  E_{\tilde\psi_n^r}}|(t(e_j))_{j\not= i}\right)\right\rbrack\;.
\end{eqnarray*}
Now, when $(t(e_j))_{j\not= i}$ is fixed, $t(e_i)\mapsto\II_{e_i\in
  E_{\tilde\psi_n^r}}$ is a non-increasing function and $t(e_i)\mapsto
t(e_i)\II_{t(e_i)\geq a}$ is of course non-decreasing. Furthermore,
since the variables $(t(e_i))$ are independent, the conditional
expectation $\EE\left(.|(t(e_j))_{j\not= i}\right)$ corresponds to
  expectation over $t(e_i)$, keeping $(t(e_j))_{j\not= i}$
  fixed. Thus, Chebyshev's association inequality (see \cite{HardyLittlewoodPolya52}, p.~43) implies:
\begin{align*}
\EE\big(t(e_i)\II_{t(e_i)\geq a}\II_{e_i\in
  E_{\tilde\psi_n^r}} & |(t(e_j))_{j\not= i}\big)\\
&\,\leq\, \EE\left(t(e_i)\II_{t(e_i)\geq a}|(t(e_j))_{j\not= i}\right)\EE\left(\II_{e_i\in
  E_{\tilde\psi_n^r}}|(t(e_j))_{j\not= i}\right)\;,\\
&\,=\, \EE\left(t(e_1)\II_{t(e_1)\geq a}\right)\EE\left(\II_{e_i\in
  E_{\tilde\psi_n^r}}|(t(e_j))_{j\not= i}\right)\;.
\end{align*}
Thus,
\begin{equation}
\label{eq:nucontinue}0\leq \EE(\psi_n^r)-\EE(\tilde\psi_n^r) \leq
\EE\left(t(e_1)\II_{t(e_1)\geq a}\right)\EE( \card (E_{\tilde\psi_n^r}))\leq
  r\EE\left(t(e_1)\II_{t(e_1)\geq a}\right)\;.
\end{equation}
Now, since $F$ has a finite moment of order 1, we can choose
$a=a(\eta,F,d)$ such that:
$$\frac{1-\eta}{\eps}\EE\left(t(e_1)\II_{t(e_1)\geq
    a}\right)\leq\frac{\eta}{2}\;,$$
to get:
\begin{eqnarray}
\nonumber
 \EE(\psi_n^r)-\EE(\tilde\psi_n^r)\leq\frac{\eta}{2}\EE(\phi_n^\kappa)\leq\frac{\eta}{2}\EE(\psi_n^r)\;,\\
\label{eq:taunrtautilde}
\PP(\psi_n^r\leq (1-\eta)\EE(\psi_n^r))\leq\PP\left(\tilde \psi_n^r\leq\EE(\tilde \psi_n^r)-\frac{\eta}{2}\EE( \psi_n^r)\right)\;.
\end{eqnarray}
Now, we shall use Corollary 3 in \cite{Boucheronetal03}. To this end, we
need some notations. We take $\tilde t'$
an independent collection of capacities with the same law as $\tilde t=(\tilde t(e_i))_{i=1\ldots,m_n}$. For each
edge $e_i\in\cyl(A,h)$, we denote by $\tilde t^{(i)}$ the collection of
capacities obtained from $\tilde t$ by replacing $\tilde t(e_i)$ by
$\tilde t'(e_i)$, and leaving
all other coordinates unchanged. Define:
$$V_-:=\EE\left\lbrack\left.\sum_{i=1}^{m_n}(\tilde\psi_n^r(t)-\tilde\psi_n^r(t^{(i)}))_-^2\right|t\right\rbrack\;,$$
where $\tilde\psi_n^r(t)$ is the maximal flow through $\cyl(nA,h(n))$ when
capacities are given by $t$. We shall denote by
$R_{\tilde\psi_n^r}$ the intersection of all the cuts whose capacity achieves the minimum in the
definition of $\tilde \psi_n^r$. Observe that:
$$\tilde\psi_n^r(t^{(i)})-\tilde\psi_n^r(t)\leq (\tilde t'(e_i)-\tilde
t(e_i))\II_{e_i\in R_{\tilde\psi_n^r}}\;,$$
and thus,
$$V_-\leq a^2 \EE [\card(R_{\tilde\psi_n^r})]\leq a^2r=a^2(1-\eta)\EE(\phi_n^\kappa)/\eps\;.$$
Thus,
Corollary 3 in \cite{Boucheronetal03} implies that, for every
$\eta\in]0,1[$,
$$\PP \left(\tilde \psi_n^r\leq \EE(\tilde \psi_n^r)-\frac{\eta}{2}\EE( \psi_n^r)\right)\leq
e^{-\frac{\EE(\psi_n^r)^2\eta^2\eps}{16a^2(1-\eta)\EE(\phi_n^\kappa)}}\leq e^{-\frac{\EE(\phi_n^\kappa)\eta^2\eps}{16a^2(1-\eta)}}\;.$$
Using inequalities (\ref{eq:taunrtautilde}) and (\ref{eq:tautaunr})
and taking the maximum over $\kappa\in D_n$, this ends the proof of Inequality (\ref{chapitre6eq:deviationkappa}). 

To see that (\ref{eq:dev_phi_ordre1}) holds, notice that
$\EE(\phi_n)\leq \min_{\kappa\in D_n}\EE(\phi_n^\kappa)$. Thus,
(\ref{eq:dev_phi_ordre1}) is a consequence of inequalities
(\ref{chapitre6eq:deviationkappa}) and (\ref{chapitre6eq:Oh2}).
\end{dem}


\section{Asymptotic behaviour of the expectation of the maximal flow}
\label{sec:CVesp}

\subsection{Upper bound}
\label{chapitre6subsec:upperbound}

From now on, we suppose that the conditions
(\ref{chapitre6eq:conditionsFLGN}) on $F$ and
(\ref{chapitre6eq:conditionshn})  on $h$ are satisfied. We consider a line segment $A$, of
orthogonal unit vector $\va = (\cos\theta, \sin\theta)$ for $\theta
\in [0, \pi[$, and a function $h: \NN \rightarrow \RR^+ $ satisfying $\lim_{n\rightarrow
  \infty}h(n) = +\infty$. Recall that $\D_n=\D(nA,h(n))$. For all $\wt \in
\D_n $,  we define
$$  k_n \,=\, \frac{1}{2} + \frac{nl(A) \tan (\wt - \theta)}{4 h(n)} \,,$$
and thus $\kappa_n=(k_n,\widetilde{\theta})\in D_n$. We want to
compare $\phi_n^{\kappa_n}$ with the maximal flow $\tau$ in a cylinder
inside $\cyl(nA, h(n))$ and oriented towards the direction
$\widetilde{\theta}$. In fact, we must use the subadditivity of $\tau$ and
compare $\phi_n^{\kappa_n}$ with a sum of such variables $\tau$.

\begin{figure}[!ht]
\centering
\begin{picture}(0,0)%
\includegraphics{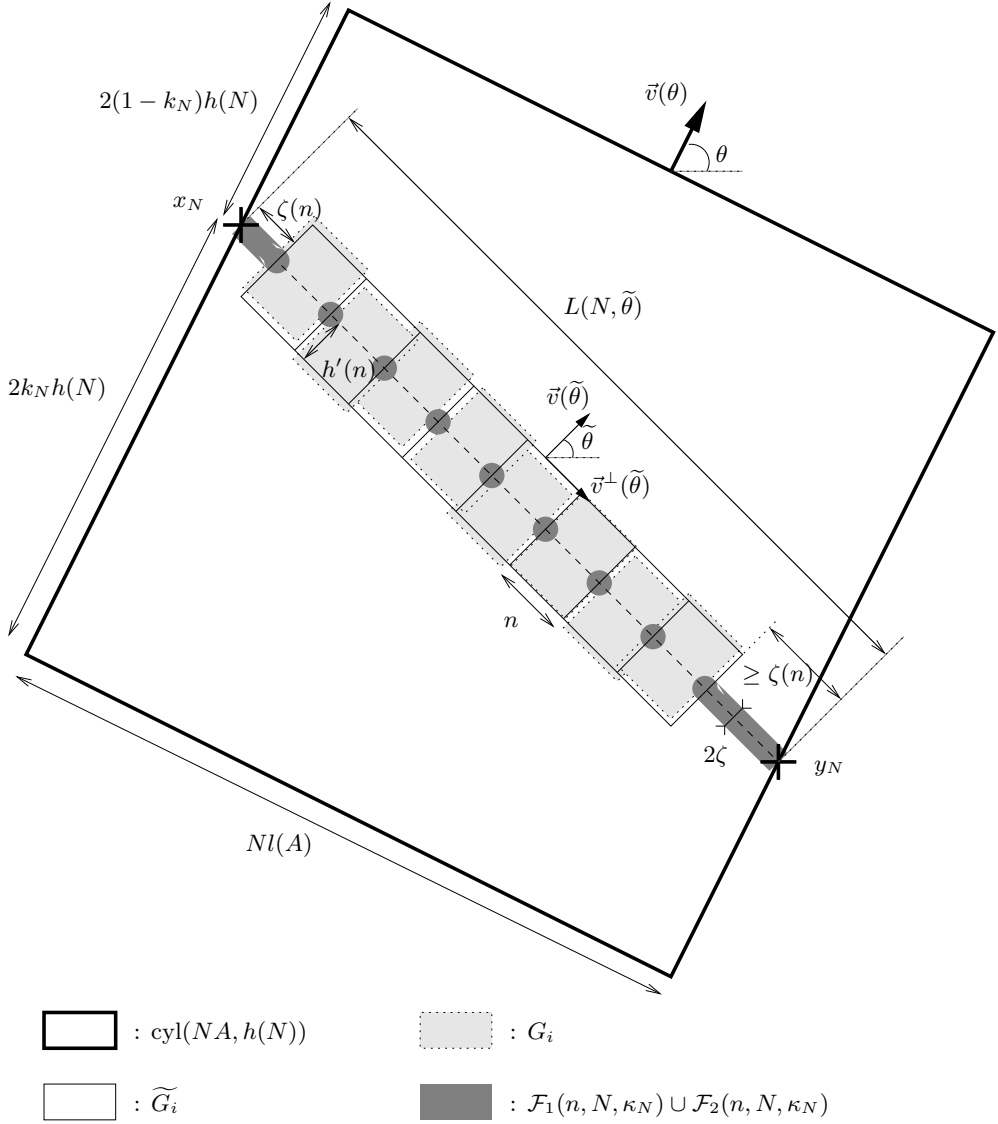}%
\end{picture}%
\setlength{\unitlength}{2960sp}%
\begingroup\makeatletter\ifx\SetFigFont\undefined%
\gdef\SetFigFont#1#2#3#4#5{%
  \reset@font\fontsize{#1}{#2pt}%
  \fontfamily{#3}\fontseries{#4}\fontshape{#5}%
  \selectfont}%
\fi\endgroup%
\begin{picture}(8295,9399)(1039,-10573)
\put(3301,-2986){\makebox(0,0)[lb]{\smash{{\SetFigFont{9}{10.8}{\rmdefault}{\mddefault}{\updefault}{\color[rgb]{0,0,0}$\zeta(n)$}%
}}}}
\put(1876,-4486){\makebox(0,0)[rb]{\smash{{\SetFigFont{9}{10.8}{\rmdefault}{\mddefault}{\updefault}{\color[rgb]{0,0,0}$2k_Nh(N)$}%
}}}}
\put(3151,-2086){\makebox(0,0)[rb]{\smash{{\SetFigFont{9}{10.8}{\rmdefault}{\mddefault}{\updefault}{\color[rgb]{0,0,0}$2(1-k_N)h(N)$}%
}}}}
\put(2701,-2911){\makebox(0,0)[rb]{\smash{{\SetFigFont{9}{10.8}{\rmdefault}{\mddefault}{\updefault}{\color[rgb]{0,0,0}$x_N$}%
}}}}
\put(6976,-2536){\makebox(0,0)[lb]{\smash{{\SetFigFont{9}{10.8}{\rmdefault}{\mddefault}{\updefault}{\color[rgb]{0,0,0}$\theta$}%
}}}}
\put(6751,-2011){\makebox(0,0)[rb]{\smash{{\SetFigFont{9}{10.8}{\rmdefault}{\mddefault}{\updefault}{\color[rgb]{0,0,0}$\vec{v}(\theta)$}%
}}}}
\put(7801,-7636){\makebox(0,0)[lb]{\smash{{\SetFigFont{9}{10.8}{\rmdefault}{\mddefault}{\updefault}{\color[rgb]{0,0,0}$y_N$}%
}}}}
\put(6976,-7561){\makebox(0,0)[b]{\smash{{\SetFigFont{9}{10.8}{\rmdefault}{\mddefault}{\updefault}{\color[rgb]{0,0,0}$2\zeta$}%
}}}}
\put(5851,-4936){\makebox(0,0)[lb]{\smash{{\SetFigFont{9}{10.8}{\rmdefault}{\mddefault}{\updefault}{\color[rgb]{0,0,0}$\widetilde{\theta}$}%
}}}}
\put(5926,-4561){\makebox(0,0)[rb]{\smash{{\SetFigFont{9}{10.8}{\rmdefault}{\mddefault}{\updefault}{\color[rgb]{0,0,0}$\vec{v}(\wt)$}%
}}}}
\put(5926,-5311){\makebox(0,0)[lb]{\smash{{\SetFigFont{9}{10.8}{\rmdefault}{\mddefault}{\updefault}{\color[rgb]{0,0,0}$\vec{v}^{\bot}(\wt)$}%
}}}}
\put(5326,-6436){\makebox(0,0)[rb]{\smash{{\SetFigFont{9}{10.8}{\rmdefault}{\mddefault}{\updefault}{\color[rgb]{0,0,0}$n$}%
}}}}
\put(7801,-6886){\makebox(0,0)[rb]{\smash{{\SetFigFont{9}{10.8}{\rmdefault}{\mddefault}{\updefault}{\color[rgb]{0,0,0}$\geq \zeta(n)$}%
}}}}
\put(5701,-3811){\makebox(0,0)[lb]{\smash{{\SetFigFont{9}{10.8}{\rmdefault}{\mddefault}{\updefault}{\color[rgb]{0,0,0}$L(N,\widetilde{\theta})$}%
}}}}
\put(5251,-10486){\makebox(0,0)[lb]{\smash{{\SetFigFont{9}{10.8}{\rmdefault}{\mddefault}{\updefault}{\color[rgb]{0,0,0}: $\F_1(n,N,\kappa_N) \cup \F_2(n,N,\kappa_N)$}%
}}}}
\put(5251,-9886){\makebox(0,0)[lb]{\smash{{\SetFigFont{9}{10.8}{\rmdefault}{\mddefault}{\updefault}{\color[rgb]{0,0,0}: $G_i$}%
}}}}
\put(2101,-10486){\makebox(0,0)[lb]{\smash{{\SetFigFont{9}{10.8}{\rmdefault}{\mddefault}{\updefault}{\color[rgb]{0,0,0}: $\widetilde{G_i}$}%
}}}}
\put(2101,-9886){\makebox(0,0)[lb]{\smash{{\SetFigFont{9}{10.8}{\rmdefault}{\mddefault}{\updefault}{\color[rgb]{0,0,0}: $\cyl(NA,h(N))$}%
}}}}
\put(3601,-8311){\makebox(0,0)[rb]{\smash{{\SetFigFont{9}{10.8}{\rmdefault}{\mddefault}{\updefault}{\color[rgb]{0,0,0}$N l(A)$}%
}}}}
\put(3676,-4336){\makebox(0,0)[lb]{\smash{{\SetFigFont{9}{10.8}{\rmdefault}{\mddefault}{\updefault}{\color[rgb]{0,0,0}$h'(n)$}%
}}}}
\end{picture}%
\caption{The cylinders $\cyl(NA,h(N))$ and $G_i$, for $i=1,...,\M$.}
\label{chapitre6emboitement1}
\end{figure}

We consider $n$ and $N$ in $\NN$, with $N$ a lot bigger than $n$. The following definitions can seem a little bit
complicated, but Figure \ref{chapitre6emboitement1} is more explicit.
We choose two functions $h', \zeta :\NN \rightarrow \RR^+$ such that
$$\lim_{n\rightarrow \infty} h'(n) \,=\, \lim_{n\rightarrow \infty}
\zeta(n) \,=\, + \infty \,,$$
and 
\begin{equation}
\label{chapitre6cond1}
 \lim_{n\rightarrow \infty}\frac{h'(n)}{\zeta (n)} \,=\, 0 \,. 
\end{equation}
We consider a fixed $\wt \in \D_N$. Let 
$$\vb \,=\, (\cos \widetilde{\theta}, \sin \widetilde{\theta}) \quad
\textrm{and} \quad \vd \,=\, (\sin
\widetilde{\theta}, -\cos \widetilde{\theta}) \,.$$
In $\cyl(NA,h(N))$, we denote by $x_N$ and $y_N$ the two points
corresponding to the boundary conditions $\kappa_N$, such that
$\overrightarrow{x_Ny_N} \cdot \vd >0 $. Notice that according to our
choice of $k_N$, the segments $[x_N,y_N]$ and $NA$ cut each other in their middle. If we denote by $L(N,\wt)$ the
distance between $x_N$ and $y_N$, we have:
$$ L(N,\wt)\,=\, \frac{Nl(A)}{\cos(\widetilde{\theta}-\theta)} \,. $$
We define
$$ \cyl'(n) \,=\, \cyl ([0 ,n \vd ],h'(n))  \,.$$
We will translate $\cyl '(n)$ numerous times inside $\cyl(NA,h(N))$. We
define
$$ t_i \,=\, x_N + \left(\zeta(n) + (i-1) n \right) \vd \,,$$
for $i=1,..., \M$, where
$$ \M \,=\, \M(n,N) \,=\, \left\lfloor \frac{L(N,\wt)  -2\zeta(n)}{n} \right\rfloor \,.$$
Of course we consider only $N$ large enough to have $\M \geq 2$. For
$i=1,...,\M$, we denote by $\widetilde{G_i}$ the image of $\cyl'(n)$ by the
translation of vector $\overrightarrow{0 t_i}$. For $n$ (and thus $N$)
sufficiently large, thanks to condition (\ref{chapitre6cond1}), we know that
$\widetilde{G_i} \subset \cyl(NA,h(N))$ for all $i$. We can translate $\widetilde{G_i}$ again by a vector of
norm strictly smaller than $1$ to obtain an integer translate of $\cyl'(n)$
(i.e., a translate by a vector whose coordinates are in $\ZZ^2$) that we
will call $G_i$. Now we want to glue together cutsets of boundary condition
$(1/2, \wt)$ in the cylinders $G_i$. We define:
$$ \F_1(n,N,\kappa_N) \,=\, \left(\bigcup_{i=1}^{\M} \mathcal{V} (t_i,
  \zeta_0) \right) \, \bigcap \,\cyl(NA,h(N))  \,,$$
where $\zeta_0$ is a fixed constant larger than $4$, and:
$$ \F_2(n,N,\kappa_N) \,=\, \mathcal{V} \left( [x_N,x_N + \zeta(n) \vd ]
  \cup [z_\M , y_N ] ,\zeta_0 \right)\, \bigcap \,\cyl(NA,h(N)) \,.$$
Let $F_1(n,N,\kappa_N)$ (respectively $F_2(n,N,\kappa_N)$) be the set of
the edges included in $\F_1(n,N,\kappa_N)$ (respectively
$\F_2(n,N,\kappa_N)$). If for every $i=1,...,\M$, $\G_i$ is a cutset of
boundary condition $(1/2,\wt)$ in $G_i$, then
$$ \bigcup_{i=1}^{\M} \G_i \cup F_1(n,N,\kappa_N) \cup F_2(n,N,\kappa_N) $$
contains a cutset of boundary conditions $\kappa_N$ in $\cyl(NA,h(N))$. We
obtain:
\begin{equation}
\label{chapitre6lien_phiF_tau_1} 
\phi_N^{\kappa_N} \,\leq\, \sum_{i=1}^{\M} \tau(G_i,\vb ) +
V(F_1(n,N,\kappa_N) \cup F_2(n,N,\kappa_N) ) \,,
\end{equation}
and so,
\begin{equation}
\label{chapitre6lien_phi_tau} 
\forall \widetilde{\theta}\in \D_N \qquad \phi_N \,\leq\,
\phi_{N}^{\kappa_N} \,\leq\, \sum_{i=1}^{\M}\tau(G_i, \vb ) + V( F_1(n,N,\kappa_N) \cup F_2(n,N,\kappa_N) ) \,.
\end{equation}
There exists a constant $C_5$ such that:
$$ \card ( F_1(n,N,\kappa_N) ) \,\leq\, C_5
  \M  \quad \textrm{and} \quad \card (F_2(n,N,\kappa_N))\,\leq\,  C_5 \left( \zeta(n) + n \right) \,,$$
and since the set of edges $F_1(n,N,\kappa_N) \cup F_2(n,N,\kappa_N)$ is deterministic,
$$ \EE [V(F_1(n,N,\kappa_N) \cup F_2(n,N,\kappa_N))] \,\leq\, C_5 \EE(t)\left(
  \M +  \zeta(n) + n \right) \,. $$
So
\begin{equation}
\label{eqNn}
\forall \widetilde{\theta}\in \D_N \qquad
\frac{\EE(\phi_N)}{Nl(A)} \,\leq\, \frac{\M n}{N l(A)}
\times \frac{\EE [\tau(\cyl'(n), \vb )]}{n}
+ \frac{C_5 \EE(t)\left( \M +  \zeta(n) + n \right) }{Nl(A)} \,.
\end{equation}
We want to send $N$ to infinity. First, let $\wt \in \underline{\D}$. Then
for all $N$ large enough, $\wt \in \D_N$, and thus for all $n$ large enough
we have
$$ \limsup_{N\rightarrow \infty}\frac{\EE(\phi_N)}{Nl(A)} \,\leq\,
\frac{1}{\cos(\wt - \theta)}  \frac{\EE [\tau(\cyl'(n), \vb )]}{n} +
\frac{C_5 \EE(t)}{n \cos(\wt - \theta)}\,.   $$
Sending $n$ to infinity, thanks to Theorem \ref{chapitre6thm:LGNtau}, we
obtain that
\begin{equation}
\label{chapitre6upper1}
 \limsup_{N\rightarrow \infty } \frac{\EE(\phi_N)}{Nl(A)}
\,\leq\, \inf_{\widetilde{\theta} \in \underline{\D}}
\frac{\nu_{\widetilde{\theta}}}{\cos (\widetilde{\theta}-\theta)} \,.
\end{equation}
We now suppose that $\wt \in \overline{\D}$. Let $\psi:\NN \rightarrow
\NN$ be strictly increasing and such that for all $N$, $\wt \in \D_{\psi(N)}$. Then
thanks to Equation (\ref{eqNn}), sending first $N$ to infinity and then $n$
to infinity, we obtain that
\begin{equation}
\label{chapitre6eq:upper1overline}
\liminf_{N\rightarrow \infty } \frac{\EE(\phi_N)}{Nl(A)}
\,\leq\, \limsup_{N\rightarrow \infty } \frac{\EE(\phi_{\psi(N)})}{\psi(N)l(A)}
\,\leq\, \inf_{\widetilde{\theta} \in \overline{\D}}
\frac{\nu_{\widetilde{\theta}}}{\cos (\widetilde{\theta}-\theta)} \,.
\end{equation}


\subsection{Lower bound}
\label{chapitre6subsec:lowerbound}

\begin{figure}[!ht]
\centering
\begin{picture}(0,0)%
\includegraphics{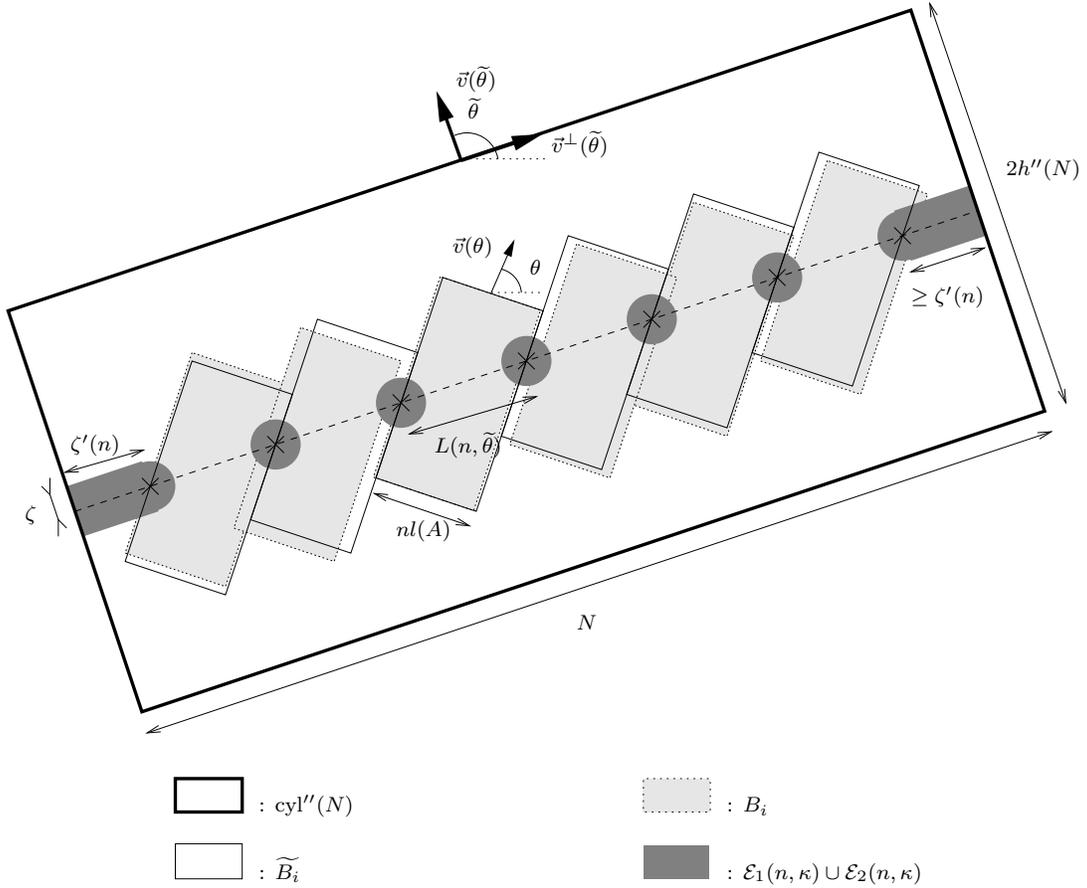}%
\end{picture}%
\setlength{\unitlength}{2763sp}%
\begingroup\makeatletter\ifx\SetFigFont\undefined%
\gdef\SetFigFont#1#2#3#4#5{%
  \reset@font\fontsize{#1}{#2pt}%
  \fontfamily{#3}\fontseries{#4}\fontshape{#5}%
  \selectfont}%
\fi\endgroup%
\begin{picture}(9540,7965)(1168,-8539)
\put(6301,-6226){\makebox(0,0)[lb]{\smash{{\SetFigFont{8}{9.6}{\rmdefault}{\mddefault}{\updefault}{\color[rgb]{0,0,0}$N$}%
}}}}
\put(10156,-2146){\makebox(0,0)[lb]{\smash{{\SetFigFont{8}{9.6}{\rmdefault}{\mddefault}{\updefault}{\color[rgb]{0,0,0}$2h''(N)$}%
}}}}
\put(1456,-5236){\makebox(0,0)[rb]{\smash{{\SetFigFont{8}{9.6}{\rmdefault}{\mddefault}{\updefault}{\color[rgb]{0,0,0}$\zeta$}%
}}}}
\put(2206,-4621){\makebox(0,0)[rb]{\smash{{\SetFigFont{8}{9.6}{\rmdefault}{\mddefault}{\updefault}{\color[rgb]{0,0,0}$\zeta'(n)$}%
}}}}
\put(3451,-8461){\makebox(0,0)[lb]{\smash{{\SetFigFont{8}{9.6}{\rmdefault}{\mddefault}{\updefault}{\color[rgb]{0,0,0}: $\widetilde{B_i}$}%
}}}}
\put(7651,-7861){\makebox(0,0)[lb]{\smash{{\SetFigFont{8}{9.6}{\rmdefault}{\mddefault}{\updefault}{\color[rgb]{0,0,0}: $B_i$}%
}}}}
\put(9301,-3286){\makebox(0,0)[lb]{\smash{{\SetFigFont{8}{9.6}{\rmdefault}{\mddefault}{\updefault}{\color[rgb]{0,0,0}$\geq \zeta'(n)$}%
}}}}
\put(7651,-8461){\makebox(0,0)[lb]{\smash{{\SetFigFont{8}{9.6}{\rmdefault}{\mddefault}{\updefault}{\color[rgb]{0,0,0}: $\E_1(n,\kappa)\cup\E_2(n,\kappa)$}%
}}}}
\put(3451,-7861){\makebox(0,0)[lb]{\smash{{\SetFigFont{8}{9.6}{\rmdefault}{\mddefault}{\updefault}{\color[rgb]{0,0,0}: $\cyl''(N)$}%
}}}}
\put(5221,-1351){\makebox(0,0)[lb]{\smash{{\SetFigFont{8}{9.6}{\rmdefault}{\mddefault}{\updefault}{\color[rgb]{0,0,0}$\vec{v}(\wt)$}%
}}}}
\put(5551,-2836){\makebox(0,0)[rb]{\smash{{\SetFigFont{8}{9.6}{\rmdefault}{\mddefault}{\updefault}{\color[rgb]{0,0,0}$\vec{v}(\theta)$}%
}}}}
\put(5326,-1636){\makebox(0,0)[lb]{\smash{{\SetFigFont{8}{9.6}{\rmdefault}{\mddefault}{\updefault}{\color[rgb]{0,0,0}$\widetilde{\theta}$}%
}}}}
\put(6076,-1936){\makebox(0,0)[lb]{\smash{{\SetFigFont{8}{9.6}{\rmdefault}{\mddefault}{\updefault}{\color[rgb]{0,0,0}$\vec{v}^{\bot}(\wt)$}%
}}}}
\put(5881,-3031){\makebox(0,0)[lb]{\smash{{\SetFigFont{8}{9.6}{\rmdefault}{\mddefault}{\updefault}{\color[rgb]{0,0,0}$\theta$}%
}}}}
\put(5026,-4636){\makebox(0,0)[lb]{\smash{{\SetFigFont{8}{9.6}{\rmdefault}{\mddefault}{\updefault}{\color[rgb]{0,0,0}$L(n,\widetilde{\theta})$}%
}}}}
\put(5176,-5386){\makebox(0,0)[rb]{\smash{{\SetFigFont{8}{9.6}{\rmdefault}{\mddefault}{\updefault}{\color[rgb]{0,0,0}$nl(A)$}%
}}}}
\end{picture}%
\caption{The cylinders $\cyl''(N)$ and $B_i$, for $i=1,...,\N$.}
\label{chapitre6emboitement2}
\end{figure}

We do the symmetric construction of the one done in section
\ref{chapitre6subsec:upperbound}. We consider $n$ and $N$ in $\NN$ and take $N$ a lot bigger than $n$. We
choose functions $\zeta', h'': \NN \rightarrow \RR^+$ such that
$$ \lim_{n\rightarrow \infty} \zeta'(n) \,=\, \lim_{n\rightarrow \infty}
h''(n) \,=\, +\infty \,, $$
and
\begin{equation}
\label{chapitre6cond2}
 \lim_{n\rightarrow \infty}\frac{ h(n)}{\zeta' (n)} \,=\,0 \,. 
\end{equation}
We consider $\kappa=(k,\widetilde{\theta})\in D_n$. Keeping the same
notations as in section \ref{chapitre6subsec:upperbound}, we define
$$ \cyl''(N) \,=\, \cyl \left([0,N \vd ] ,
  h''(N)\right) \,. $$
We will translate
$\cyl (nA, h(n))$ numerous times in $\cyl ''(N)$. The figure
\ref{chapitre6emboitement2} is more explicit than the following definitions.
The condition $\kappa$ defines two points $x_n$ and $y_n$ on the boundary of
$\cyl(nA, h(n))$ (see  section \ref{chapitre6subsec:upperbound}). As in section \ref{chapitre6subsec:upperbound}, we denote
by $L(n,\wt)$ the distance between $x_n$ and $y_n$, and we have
$$ L(n,\wt) \,=\, \frac{nl(A)}{cos(\widetilde{\theta}-\theta)}\,. $$ We
define
$$ z_i \,=\, \left(\zeta'(n) + (i-1) L(n,\wt)\right) \vd \,,$$
for $i=1,..., \N$, where
$$ \N \,=\, \left\lfloor \frac{N-2\zeta'(n)}{L(n,\wt)} \right\rfloor \,.$$
Of course we consider only $N$ large enough to have $\N \geq 2$. For
$i=1,...,\N$, we denote by $\widetilde{B_i}$ the image of $\cyl
(nA, h(n))$ by the translation of vector $\overrightarrow{x_n z_i}$. For $N$
sufficiently large, thanks to condition (\ref{chapitre6cond2}), we know that
$\widetilde{B_i} \subset \cyl''(N)$ for all $i$. We can translate $\widetilde{B_i}$ again by a vector of
norm strictly smaller than $1$ to obtain an integer translate of $\cyl(nA,
h(n))$ (i.e., a translate by a vector whose coordinates are in $\ZZ^2$) that we will call $B_i$. Now we want to glue together
cutsets of boundary condition $\kappa$ in the different $B_i$'s. We define:
$$ \E_1(n,N,\kappa) \,=\, \left(\bigcup_{i=1}^\N \mathcal{V} (z_i, \zeta) \right)
\, \bigcap \,\cyl''(N)  \,,$$
where $\zeta$ is still a fixed constant bigger than $4$, and:
$$ \E_2(n,N,\kappa) \,=\, \mathcal{V} \left( [0,\zeta'(n)
  \vd ] \cup [z_\N, N \vd ]
  ,\zeta \right)\, \bigcap \,\cyl''(N) \,.$$
Let $E_1(n,N,\kappa)$ (respectively $E_2(n,N,\kappa)$) be the set of the edges included in
$\E_1(n,N,\kappa)$ (respectively $\E_2(n,N,\kappa)$).
Then, still by gluing cutsets together, we obtain:
\begin{equation}
\label{chapitre6lien_phiF_tau_2}
\tau(\cyl''(N), \vb ) \,\leq\, \sum_{i=1}^\N
\phi^{\kappa}(B_i, \va) + V(E_1(n,N,\kappa) \cup E_2(n,N,\kappa)) \,.
\end{equation}
On one hand, there exists a constant $C_6$ (independent of $\kappa$) such that:
$$ \card (E_1(n,N,\kappa) \cup E_2(n,N,\kappa)) \,\leq\, C_6 \left( \N + \zeta'(n) +
  L(n,\wt) \right) \,,$$ 
and since the sets $E_1(n,N,\kappa)$ and $E_2(n,N,\kappa)$ are deterministic, we deduce:
$$ \EE [ V(E_1(n,N,\kappa) \cup E_2(n,N,\kappa))] \,\leq\, C_6 \EE(t) \left( \N + \zeta'(n) +
  L(n,\wt) \right)\,.  $$
On the other hand, the variables $(\phi^{\kappa}(B_i))_{i=1,...,\N}$ are
identically distributed, with the same law as $\phi_{n}^\kappa$ (because we
only consider integer translates), so (\ref{chapitre6lien_phiF_tau_2}) leads to
$$ \EE [\tau(\cyl''(N), \vb )] \,\leq\, \N \EE
[\phi_{n}^\kappa] +  C_6 \EE(t) \left( \N + \zeta'(n) + L(n,\wt) \right)\,.  $$
Dividing by $N$ and sending $N$ to infinity, we get, thanks to Theorem
\ref{chapitre6thm:LGNtau}:
$$ \nu_{\widetilde{\theta}} \,\leq\, \frac{\EE[\phi_{n}^\kappa]}{L(n,\wt)} +
\frac{C_6 \EE(t)}{ L(n,\wt)}\,, $$
and so:
$$ \frac{\EE [\phi_{n}^\kappa]}{nl(A)} \,\geq\,
\frac{\nu_{\widetilde{\theta}}}{\cos (\widetilde{\theta} - \theta)} -
\frac{C_6 \EE(t)}{nl(A)} \,.$$
Since $C_6$ is independent of $\kappa$,
$$ \inf_{\kappa\in D_n} \frac{\EE [\phi_{n}^\kappa]}{nl(A)} \,\geq\,
\inf_{\widetilde{\theta} \in \D_n} \frac{\nu_{\widetilde{\theta}}}{\cos
  (\widetilde{\theta} - \theta)} - \frac{C_6 \EE(t)}{nl(A)} \,.$$
First, we affirm:
\begin{equation}
\label{chapitre6liminf_inf}
 \liminf_{n\rightarrow \infty}
\inf_{\widetilde{\theta} \in \D_n} \frac{\nu_{\widetilde{\theta}}}{\cos
  (\widetilde{\theta} - \theta)} \,\geq\, \inf_{\widetilde{\theta} \in
  \overline{\D}}\frac{\nu_{\widetilde{\theta}}}{\cos
  (\widetilde{\theta} - \theta)}\,,
\end{equation}
and thus:
\begin{equation}
\label{chapitre6lower1}
 \liminf_{n\rightarrow \infty} \inf_{\kappa\in D_n} \frac{\EE
   [\phi_{n}^\kappa]}{nl(A)} \,\geq\,  \inf_{\widetilde{\theta} \in
  \overline{\D}}\frac{\nu_{\widetilde{\theta}}}{\cos
  (\widetilde{\theta} - \theta)}\,.
\end{equation}
We also claim that:
\begin{equation}
\label{chapitre6limsup_inf}
 \limsup_{n\rightarrow \infty}
\inf_{\widetilde{\theta} \in \D_n} \frac{\nu_{\widetilde{\theta}}}{\cos
  (\widetilde{\theta} - \theta)} \,\geq\, \inf_{\widetilde{\theta} \in
  \underline{\D}}\frac{\nu_{\widetilde{\theta}}}{\cos
  (\widetilde{\theta} - \theta)}\,,
\end{equation}
and therefore:
\begin{equation}
\label{chapitre6eq:lower1underline}
 \limsup_{n\rightarrow \infty} \inf_{\kappa\in D_n} \frac{\EE
   [\phi_{n}^\kappa]}{nl(A)} \,\geq\,  \inf_{\widetilde{\theta} \in
  \underline{\D}}\frac{\nu_{\widetilde{\theta}}}{\cos
  (\widetilde{\theta} - \theta)}\,.
\end{equation}

Let us prove  Inequality (\ref{chapitre6liminf_inf}). In fact, we will state a
more general result:
\begin{lem}
\label{chapitre6inversion}
Let $\theta \in [0,\pi[$, and $f$ be a lower semi-continuous function from
$[\theta-\pi/2, \theta+\pi/2]$ to $\RR^+ \cup \{+\infty\}$. 
Then we have
$$ \liminf_{n\rightarrow \infty} \inf_{\wt \in \D_n} f(\wt) \,\geq\,
\inf_{\wt \in \\ad (\overline{\D})} f(\wt) \,, $$
where $\ad (\overline{\D})$ is the adherence of $\overline{\D}$.
\end{lem}

\begin{dem}
We consider a positive
$\varepsilon$. For all $n$, since $f$ is lower semi-continuous and $\D_n$
is compact, there exists $\wt_n \in \D_n$ such that
$f (\wt_n) = \inf_{\wt \in \D_n} f(\wt) $. Up to
extracting a subsequence, we can suppose that the sequence $(\inf_{\wt \in
  \D_n}f (\wt))_{n\geq 0}$
converges towards $\displaystyle{\liminf_{n\rightarrow \infty} \inf_{\wt\in \D_n}
f(\wt)}$, and so:
$$\lim_{n\rightarrow \infty} f( \wt_n ) \,=\,
\liminf_{n\rightarrow \infty} \inf_{\wt\in \D_n} f(\wt) \,.$$
The
sequence $(\wt_n)_{n\geq 0}$ (in fact the previous subsequence) takes values in the
compact  $[\theta-\pi/2, \theta + \pi/2]$, so up to extracting a second
subsequence we can suppose that $(\wt_n)_{n\geq 0}$ converges towards a limit
$\wt_\infty$ in this compact. Since $f$ is lower semi-continuous,
$$ f(\wt_\infty) \,\leq\, \lim_{n\rightarrow \infty}f (\wt_n) \,=\,
\liminf_{n\rightarrow \infty} \inf_{\wt\in \D_n} f(\wt) \,,$$
and we just have to prove that $\wt_\infty$ belongs to
$\ad(\overline{\D})$. Indeed, for all positive $\varepsilon$, $\wt_n
\in [\wt_\infty - \varepsilon, \wt_\infty + \varepsilon ]$ for an infinite
number of $n$. We remember that all the $\D_n$ are closed intervals
centered at $\theta$. If $\wt_\infty = \theta$, the result is obvious. We
suppose that $\wt_\infty >\theta$ for example, and thus, for $\varepsilon$
small enough, $\wt_\infty - \varepsilon >\theta$. Then $[\theta, \wt_\infty -
\varepsilon]$ is included in an infinite number of $\D_n$, so $\wt_\infty -
\varepsilon$ belongs to $\overline{\D}$, and then $\wt_\infty$ belongs to
$\ad (\overline{\D})$. The same holds if $\wt_\infty < \theta$. This ends
the proof of Lemma \ref{chapitre6inversion}.
\end{dem}

We use Lemma \ref{chapitre6inversion} with $f(\wt) =
\nu_{\widetilde{\theta}} / \cos (\widetilde{\theta} - \theta)$. Here $f$ is
lower semi-continuous, because $\wt \rightarrow \nu_{\wt}$ is continuous
since it satisfies the weak triangle inequality. Indeed, it is obvious in
dimension $2$ because $\nu_{\wt} = \mu(\vb)$ which satisfies the (ordinary)
triangle inequality, but it has also been proved in any dimension $d\geq 2$
(see section 4.4 in \cite{RossignolTheret08b}). Moreover we know that $f$ is finite and continuous
  on $]\theta - \pi/2, \theta + \pi/2[$, infinite at $\theta + \pi/2$ and
  $\theta - \pi/2$ and 
$$\lim_{\wt \rightarrow \theta+\pi/2} f(\wt) \,=\,
  \lim_{\wt \rightarrow \theta-\pi/2} f(\wt) \,=\, +\infty \,,$$
so we can even say in this case:
$$\inf_{\wt \in \ad(\overline{\D})} f(\wt) = \inf_{\wt \in \overline{\D}}
f(\wt) \,,$$
and we obtain Inequality (\ref{chapitre6liminf_inf}).

Let us now prove Inequality (\ref{chapitre6limsup_inf}). We state again a more
general result:
\begin{lem}
\label{chapitre6inversion-bis}
Let $\theta \in [0,\pi[$, and $f$ be a lower semi-continuous function from
$[\theta - \pi/2, \theta+\pi/2]$ to $\RR^+ \cup \{+\infty\}$. Then we have
$$ \limsup_{n\rightarrow \infty} \inf_{\wt \in \D_n} f(\wt) \,\geq\,
\inf_{\wt \in \ad (\underline{\D})} f(\wt) \,, $$
where $\ad (\underline{\D}) $ is the adherence of $\underline{\D}$.
\end{lem}

\begin{dem}
We denote $\ad(\underline{\D})$ by $[\theta -\alpha, \theta + \alpha]$. For all integer $p\geq 1$, there exists $n_p\geq
n_{p-1}$ ($n_0 = 1$) such that:
$$\theta + \alpha + 1/p \,\notin\, \D_{n_p} \qquad \mbox{and} \qquad \theta - \alpha -
1/p \,\notin\, \D_{n_p} \,,$$
thus
$$ \D_{n_p} \,\subset \, ]\theta - \alpha - 1/p , \theta + \alpha + 1/p[
\,, $$
then
\begin{align*}
\limsup_{n\rightarrow \infty} \inf_{\wt \in \D_n} f(\wt) & \,\geq\,
\limsup_{p\rightarrow \infty} \inf_{\wt \in \D_{n_p}} f(\wt) \\
& \,\geq\, \limsup_{p\rightarrow \infty} \inf_{\wt \in [\theta - \alpha -
  1/p, \theta + \alpha + 1/p]} f(\wt) \,.
\end{align*}
The function $f$ is lower semi-continuous and $[\theta - \alpha - 1/p, \theta + \alpha +
1/p]$ is compact, so for all integers $p$ there exists $\wt_p \in [\theta - \alpha -
  1/p, \theta + \alpha + 1/p]$ such that $f(\wt_p) = \inf_{\wt \in [\theta
    - \alpha - 1/p, \theta + \alpha + 1/p]} f(\wt)$. Up to extraction, we
  can suppose that $(\wt_p)_{p\geq 1}$ converges towards a limit
  $\wt_\infty$, that belongs obviously to $[\theta-\alpha, \theta +
  \alpha]$. Finally, because $f$ is lower semi-continuous,
$$ \inf_{\wt \in [\theta-\alpha, \theta + \alpha]} f(\wt) \,\leq\,
f(\wt_\infty) \,\leq\, \limsup_{p\rightarrow \infty} f(\wt_p) \,\leq\,
\limsup_{n\rightarrow \infty} \inf_{\wt \in \D_n} f(\wt) \,,   $$
so Lemma \ref{chapitre6inversion-bis} is proved.
\end{dem}

As previously, we use Lemma \ref{chapitre6inversion-bis} with $f(\wt) =
\nu_{\widetilde{\theta}} / \cos (\widetilde{\theta} - \theta)$. Again,
we have:
$$\inf_{\wt \in \ad(\underline{\D})} f(\wt) = \inf_{\wt \in \underline{\D}}
f(\wt) \,,$$
and Equation (\ref{chapitre6limsup_inf}) is proved.


\subsection{End of the study of the mean}
\label{chapitre6subsec:endLGN}

Now, we are able to conclude the proof of (\ref{eq:liminfmoyenne}) and
(\ref{eq:limsupmoyenne}). First, we show that $\EE(\phi_n)$ and
$\min_\kappa\EE(\phi_n^\kappa)$ are of the same order.
\begin{lem}
\label{chapitre6lem:infmoyenne}
Let $A$ be a line segment in $\RR^2$. Suppose that conditions (\ref{chapitre6eq:conditionshn})
and (\ref{chapitre6eq:conditionsFLGN}) are satisfied. Then,
$$\lim_{n\rightarrow\infty}\frac{\EE(\phi_{n})}{\min_{\kappa\in
    D_n}\EE(\phi_{n}^\kappa)}=1\;.$$
\end{lem}
\begin{dem}
Notice that $\EE(\phi_{n})\leq \min_{\kappa\in
    D_n}\EE(\phi_{n}^\kappa)$, and thus it is sufficient to show that:
$$\liminf_{n\rightarrow\infty}\frac{\EE(\phi_{n})}{\min_{\kappa\in
    D_n}\EE(\phi_{n}^\kappa)}\geq 1\;.$$
Recall from (\ref{chapitre6eq:Oh2}) and Lemma \ref{chapitre6lem:duality} that there is a finite subset $\tilde D_n$ of $D_n$, such that:
$$\card(\tilde D_n) \leq C_4 h(n)^2\;,$$
for some constant $C_4$ and every $n$, and
\begin{equation}
\label{chapitre6eqphitau}\phi_{n}= \min_{\kappa\in \tilde D_n}\phi^\kappa_{n} \;.
\end{equation}
Thus, for $\eta$ in $]0,1[$,
\begin{eqnarray*}
\PP(\min_{\kappa\in \tilde D_n}\phi^\kappa_n\geq \min_{\kappa\in \tilde
  D_n}\EE(\phi^\kappa_n)(1-\eta))&=&1-\PP(\exists \kappa\in \tilde D_n,\;\phi^\kappa_n<
\min_{\kappa\in \tilde D_n}\EE(\phi^\kappa_n)(1-\eta))\;,\\
&\geq &1- |\tilde D_n|\max_{\kappa\in \tilde D_n}\PP(\phi^\kappa_n <
\min_{\kappa\in \tilde D_n}\EE(\phi^\kappa_n)(1-\eta))\;,\\
&\geq &1- C_4 h(n)^2\max_{\kappa\in \tilde D_n}\PP(\phi^\kappa_n < \EE(\phi^\kappa_n)(1-\eta))\;.
\end{eqnarray*}
Now, Proposition \ref{prop:deviation} implies that for $\eta$ in $]0,1[$,
\begin{eqnarray*}
\PP(\min_{\kappa\in D_n}\phi^\kappa_n\geq \min_{\kappa\in
  D_n}\EE(\phi^\kappa_n)(1-\eta))&\geq &1- C_4
K_1h(n)^2e^{-C(\eta,F)\min_{\kappa\in  D_n}\EE(\phi^\kappa_n)}\;,
\end{eqnarray*}
where $C(\eta,F)$ is strictly positive. Now, let $\eta_0$ be fixed in  $]0,1/2[$.
\begin{eqnarray*}
\EE(\min_{\kappa\in D_n}\phi_n^\kappa)&=&\int_0^{+\infty}\PP(\min_{\kappa\in D_n}\phi_n^\kappa
\geq t)\;dt\;,\\ 
&\geq &\int_0^{\min_{\kappa\in D_n}\EE(\phi_n^\kappa)}\PP\left(\min_{\kappa\in D_n}\phi_n^\kappa \geq \min_{\kappa\in D_n}\EE(\phi_n^\kappa)-u\right)\;du\;,\\
&\geq & \min_{\kappa\in  D_n}\EE(\phi_n^\kappa)\int_{\eta_0}^{(1-\eta_0)}\PP\left(\min_{\kappa\in
    D_n}\phi_n^\kappa\geq \min_{\kappa\in  D_n}\EE(\phi_n^\kappa)(1-\eta)\right)\;d\eta\;,\\
&\geq & \min_{\kappa\in D_n}\EE(\phi_n^\kappa)(1-2\eta_0)\left(1-C_4
K_1h(n)^2e^{-C(1-\eta_0,F)\min_{\kappa\in  D_n}\EE(\phi^\kappa_n)}\right)\;.
\end{eqnarray*}
Thanks to Inequality (\ref{chapitre6lower1}), we know that there is a strictly positive
constant $C(A)$ such that:
$$\liminf_{n\rightarrow\infty}\frac{\min_{\kappa \in \tilde
  D_n }\EE(\phi_n^\kappa)}{n}\geq C(A)\;.$$
Thus, using assumption (\ref{chapitre6eq:conditionshn}), namely the
fact that $\log h(n)$ is small
compared to $n$,
$$\liminf_{n\rightarrow\infty}\frac{\EE(\phi_{n})}{\min_{\kappa\in
    D_n}\EE(\phi_{n}^\kappa)}\geq 1-2\eta_0\;.$$
Since this is true for any $\eta_0\in ]0,1/2[$, this
finishes the proof of Lemma \ref{chapitre6lem:infmoyenne}.
\end{dem}
Now, inequalities (\ref{chapitre6upper1}),
(\ref{chapitre6eq:lower1underline}) and Lemma
\ref{chapitre6lem:infmoyenne} give:
\begin{equation}
\label{cclesp1}
 \limsup_{n\rightarrow \infty} \frac{\EE [\phi_n] }{nl(A)} \,=\,
\inf_{\wt \in \underline{\D}} \frac{\nu_{\wt}}{ \cos(\wt - \theta)}
\end{equation}
which is (\ref{eq:limsupmoyenne}). Similarly, inequalities (\ref{chapitre6eq:upper1overline}),
(\ref{chapitre6lower1}) and Lemma
\ref{chapitre6lem:infmoyenne} give:
\begin{equation}
\label{cclesp2}
\liminf_{n\rightarrow \infty} \frac{\EE [\phi_n] }{nl(A)} \,=\,
\inf_{\wt \in \overline{\D}} \frac{\nu_{\wt}}{ \cos(\wt - \theta)} \,.
\end{equation}
which is (\ref{eq:liminfmoyenne}).


\section{Proof of the law of large numbers}
\label{secLLN}

Using Borel-Cantelli's Lemma and Proposition
\ref{prop:deviation}, we obtain that
$$ \liminf_{n\rightarrow \infty} \frac{\phi_n - \EE[\phi_n]}{nl(A)}
\,\geq\,0 \,, $$
and thus, using Equations (\ref{cclesp1}) and (\ref{cclesp2}), that
\begin{equation}
\label{ajout1}
\liminf_{n\rightarrow \infty} \frac{\phi_n}{nl(A)} \,\geq\, \liminf_{n\rightarrow \infty} \frac{\EE [\phi_n] }{nl(A)} \,=\,
\inf_{\wt \in \overline{\D}} \frac{\nu_{\wt}}{ \cos(\wt - \theta)}
\end{equation}
and
\begin{equation}
\label{ajout2}
 \limsup_{n\rightarrow \infty} \frac{\phi_n}{nl(A)} \,\geq\, \limsup_{n\rightarrow \infty} \frac{\EE [\phi_n] }{nl(A)} \,=\,
\inf_{\wt \in \underline{\D}} \frac{\nu_{\wt}}{ \cos(\wt - \theta)}\,.
\end{equation}
It can seem a bit strange to bound $\displaystyle{\limsup_{n\rightarrow
    \infty} \phi_n / (nl(A))}$ from below in the study of the asymptotic
behavior of $\phi_n$. The reason is the following: we do not only want to
prove the convergence of the rescaled flow $\phi_n$ in some cases, we want to
obtain a necessary and sufficient condition for this convergence to
hold. Thus we need to know exactly the values of
$\displaystyle{\limsup_{n\rightarrow \infty} \phi_n / (nl(A))}$ and
$\displaystyle{\liminf_{n\rightarrow \infty} \phi_n / (nl(A))}$. We will
prove the converse of Inequalities (\ref{ajout1}) and
(\ref{ajout2}). For that purpose we use again the geometrical
construction performed in section \ref{chapitre6subsec:upperbound}. Suppose
only for the moment that
$$ \int_{[0,+\infty[} x \,dF(x) \,<\,\infty\,. $$
Let $\wt_1 \in \underline{\D}$ be such that
$$ \frac{\nu_{\wt_1}}{\cos (\wt_1 - \theta)} \,=\, \inf_{\wt \in
  \underline{\D}} \frac{\nu_{\wt}}{\cos(\wt - \theta)} \,. $$
Such a $\wt_1$ exists, since 
$$  \inf_{\wt \in \underline{\D}} \frac{\nu_{\wt}}{\cos(\wt - \theta)}
\,=\,  \inf_{\wt \in \ad (\underline{\D}) } \frac{\nu_{\wt}}{\cos(\wt - \theta)} $$
as stated in section \ref{chapitre6subsec:lowerbound}, $ \ad
(\underline{\D})$ is compact and the function $\wt \mapsto \nu_{\wt} / \cos
(\wt - \theta)$ is lower semi-continuous. For all $N$ large enough, $\wt_1
\in \D_N$, and we only consider such large $N$. First suppose that $0$, the
origin of the graph, is the middle of $A$. Then consider $\kappa_N = (k_N,
\wt_1)$ as defined in section \ref{chapitre6subsec:upperbound}. We performed the geometrical construction of
section \ref{chapitre6subsec:upperbound}: we consider several integer translates $G_i$, for $i=1,...,\M(n,N)$,
of $\cyl'(n)$ inside $\cyl(NA, h(N))$. Since $0$ belongs to $[x_N, y_N]$,
we can construct the cylinders $G_i$ and the sets of edges
$F_1(n,N,\kappa_N)$ and $F_2(n,N,\kappa_N)$ in such a way that
$$ \forall N_1 \leq N_2 \quad (G_i)_{i=1,...,\M(n,N_1)} \,\subset\, (G_i)_{i=1,..., \M(n,N_2)} \quad \textrm{and}\quad  F_1(n,N_1,\kappa_{N_1}) \,\subset\,
F_1(n,N_2,\kappa_{N_2})  \,.$$
We use again Inequality (\ref{chapitre6lien_phi_tau}) to obtain that:
\begin{equation}
\label{eqfin}
 \frac{\phi_N}{Nl(A)} \,\leq\, \frac{n \M}{Nl(A)} \frac{1}{\M} \sum_{i=1}^{\M}
\frac{\tau(G_i,\vec{v} (\wt_1))}{n} + \frac{V(F_1(n,N,\kappa_N) )}{Nl(A)} +
\frac{V( F_2(n,N,\kappa_N))}{Nl(A)} \,.
\end{equation}
The variables $ (\tau(G_i,\wt_1),i=1,...,\M(n,N))$ are not
independent. However, each cylinder $G_i$ can intersect at most the
two other cylinders that are its neighbours, thus we can divide the family
$ (\tau(G_i,\wt_1),i=1,...,\M(n,N))$ into two families $
(\tau(G_{i},\wt_1),i\in\{1,...,\M(n,N)\} \cap P_j)$ for $j=1,2$, $P_1
= 2\NN$ and $P_2 = 2\NN +1$, such that for each $j\in \{1,2\}$, the family $
(\tau(G_{i},\wt_1),i\in\{1,...,\M(n,N)\} \cap P_j)$ is i.i.d. Since
$$\int_{[0,+\infty[} x \,dF(x) \,<\,\infty \,,$$
it is easy to see that the variable $\tau(\cyl '(n), \wt_1)$ is integrable
(we can compare this variable with the capacity of a deterministic cutset),
and we can apply the strong law of large numbers to each of the two
families of variables described above to finally obtain that
\begin{equation}
\label{eqfin1}
\lim_{N\rightarrow \infty}  \frac{n \M}{Nl(A)} \frac{1}{\M} \sum_{i=1}^{\M}
\frac{\tau(G_i,\vec{v} (\wt_1))}{n} \,=\, \frac{1}{\cos(\wt_1 - \theta)} \frac{\EE[\tau(\cyl'(n) ,
  \wt_1)]}{n} \qquad \textrm{a.s.}
\end{equation}
Up to increasing a little the sets $F_1(n,N,\kappa_N)$, we can suppose that
for all $N$, we have
$$ \card(F_1(n,N,\kappa_N)) \,=\, C_5 \M(n,N)\,, $$
and thus, by the strong law of large numbers, we obtain that
\begin{equation}
\label{eqfin2}
\lim_{N \rightarrow \infty} \frac{V(F_1(n,N,\kappa_N) )}{Nl(A)} \,=\,
\frac{C_5 \EE[t(e)]}{n \cos(\wt_1 - \theta)}\qquad \textrm{a.s.}
\end{equation}
Moreover, we know that
$$ \card(F_2 (n,N,\kappa_N)) \,\leq\, C_5 (n+\zeta(n))\,, $$
thus for all $\eta > 0$ we have
\begin{align*}
\sum_{N\in \NN^*} \PP[V(F_2(n,N,\kappa_N))\geq \eta N l(A)]
& \,\leq\,\sum_{N\in \NN^*} \PP \left[ \sum_{i=1}^{C_5 (n+\zeta(n))} t_i \geq
  \eta N l(A) \right] \\
&\,\leq\, \EE\left[1+\frac{1}{\eta l(A)}\sum_{i=1}^{C_5 (n+\zeta(n))} t_i\right]<\infty\;.
\end{align*}
where $(t_i, i\in \NN)$ is a family of i.i.d. variables with distribution
function $F$. By a simple Borel-Cantelli's Lemma, we conclude that
\begin{equation}
\label{eqfin3}
 \lim_{N \rightarrow \infty} \frac{V(F_2(n,N,\kappa_N) )}{Nl(A)} \,=\,
0 \qquad \textrm{a.s.}
\end{equation}
Combining Equations (\ref{eqfin}), (\ref{eqfin1}), (\ref{eqfin2}) and
(\ref{eqfin3}), and sending $n$ to infinity, thanks to Theorem
\ref{chapitre6thm:LGNtau} we obtain that
$$ \limsup_{N\rightarrow \infty} \frac{\phi_N}{Nl(A)} \,\leq\,
\frac{\nu_{\wt_1}}{\cos(\wt_1 - \theta)} \,=\, \inf_{\wt\in
  \underline{\D}} \frac{\nu_{\wt}}{ \cos(\wt - \theta)}\qquad \textrm{a.s.} $$
Similarly, we can choose $\wt_2 \in \overline{\D}$ satisfying
$$ \frac{\nu_{\wt_2}}{\cos(\wt_2 - \theta)} \,=\,\inf_{\wt \in
  \overline{\D}} \frac{\nu_{\wt}}{ \cos(\wt - \theta)} \,. $$
We consider a subsequence $(\psi(N), N\in \NN)$ of $\NN$ such that for all
$N$, $\wt_2 \in \D_{\psi(N)}$. If $0$ is the middle of $A$, for every $N$
we consider $k_{\psi(N)}$ as defined in section \ref{chapitre6subsec:upperbound}, and which is such
that $\kappa_{\psi(N)} =(k_{\psi(N)} , \wt_2) \in D_{\psi(N)}$ and $0$
belongs to the segments $[x_N, y_N]$ determined by the boundary condition
$\kappa_N$. Then we obtain exactly by the same methods that
$$ \liminf_{N\rightarrow \infty} \frac{\phi_N}{Nl(A)}\,\leq\,\limsup_{N\rightarrow \infty} \frac{\phi_{\psi(N)}}{\psi(N)l(A)}  \,\leq\,
\frac{\nu_{\wt_2}}{\cos(\wt_2 - \theta)} \,=\, \inf_{\wt\in
  \overline{\D}} \frac{\nu_{\wt}}{ \cos(\wt - \theta)}\qquad \textrm{a.s.} $$
If the condition on the origin $0$ of the graph is not satisfied, we
suppose that
$$ \int_{[0,+\infty[} x^2 \,dF(x) \,<\,\infty \,. $$
To obtain Equations (\ref{eqfin1}) and (\ref{eqfin2}) in the case where $0$
is the middle of $A$, we have used the strong law of large numbers. If $0$
is not the middle of $A$ we may not construct the cylinders $(G_i,i\in
\{1,...,\M(n,N)\})$ such that the same $G_i$'s appear for different $N$. Thus
we obtain cylinders $(G_i(N), i\in \{1,...,\M(n,N)\})$ that depend on
$N$. The sets $(\tau(G_i(N),
\wt_1), i\in \{ 1,...,\M(n,N) \}\cap P_j)$ (resp. $(t(e), e\in
F_1(n,N,\kappa_N))$) are families of i.i.d. random variables for a given
$N$, and $\tau(G_i(N))$ (resp. $t(e)$) has the same law whatever the value
of $i$ and $N$ (resp. whatever $e$ and $N$), but we are not in the
conditions of application of the strong law of large numbers: we consider
the behavior of a sequence of the form
$$ \left( \frac{\sum_{i=1}^{n} X_i^{(n)}}{n}\,,\,\, n\in \NN \right) \,,$$
where $(X_i^{(j)})_{i,j}$ is an array of i.d. random variables such that for each
$n$, the variables $(X_1^{(n)},...,X_n^{(n)})$ are independent. Thanks to
Theorem~3 in \cite{HsuRobbins47}, we know that such a sequence converges
a.s. towards $\EE(X_1^{(1)})$ as soon as $\EE[(X_1^{(1)})^2] <\infty$. This
theorem is based on a result of complete convergence (see Theorem~1 in
\cite{HsuRobbins47}) and a Borel-Cantelli's Lemma. If $t(e)$ admits a
moment of order $2$, the same holds for $\tau(G_i(N),\wt_1)$, thus we can
use Theorem~3 in \cite{HsuRobbins47} to get Equations (\ref{eqfin1}) and
(\ref{eqfin2}) again.
This ends the proof of Theorem \ref{chapitre6thm:lgn}.

Obviously, the condition
\begin{equation}
\label{chapitre6CS}
 \inf_{\widetilde{\theta} \in \underline{\D}} \frac{\nu_{\widetilde{\theta}}}{\cos
  (\widetilde{\theta} - \theta)} \,=\,  \inf_{\widetilde{\theta} \in
  \overline{\D}} \frac{\nu_{\widetilde{\theta}}}{\cos (\widetilde{\theta} -
  \theta)} \,:=\, \eta_{\theta,h} \,, 
\end{equation}
necessary and sufficient for the convergence a.s. of $(\phi_n /
(nl(A)))_{n\geq 0}$, is closely linked to the asymptotic
behaviour of $h(n)/n$. Indeed we know that
$$ \D_n \,=\, [\theta - \alpha_n, \theta +\alpha_n] \,,$$
where $\alpha_n = \arctan \left( \frac{2h(n)}{nl(A)} \right)$. If
$\lim_{n\rightarrow \infty} 2h(n)/(nl(A))$ exists in $\RR^+ \cup \{+\infty \}$, and
we denote it by $\tan \alpha$ ($\alpha \in [0,\pi/2]$), then $\underline{\D}$ and $\overline{\D}$
are equal to $[\theta - \alpha, \theta + \alpha]$ or $]\theta - \alpha,
\theta + \alpha[$, and we obtain that
$\eta_{\theta,h}$ exists and 
$$\eta_{\theta,h} \,=\, \inf_{\widetilde{\theta} \in [\theta - \alpha,
  \theta + \alpha]} \frac{\nu_{\widetilde{\theta}}}{\cos (\widetilde{\theta} -
  \theta)}  \,.$$
As previously, we do not care keeping $\theta + \alpha$ and $\theta -
\alpha$ in the infimum.
Then we obtain the a.s. convergence appearing in Corollary
\ref{chapitre6corollaire}. Obviously, if there exists a $\widetilde{\theta_0}$ such that 
$$  \frac{\nu_{\widetilde{\theta_0}}}{\cos (\widetilde{\theta_0} -
  \theta)} \,=\, \inf_{\widetilde{\theta} \in [\theta-\pi/2,\theta +\pi/2]}
\frac{\nu_{\widetilde{\theta}}}{\cos (\widetilde{\theta} - \theta)} $$
and if
$$ \liminf_{n\rightarrow \infty} \frac{2h(n)}{nl(A)} \,\geq\, |\tan
(\widetilde{\theta_0} -\theta)| \,,$$
then $\eta_{\theta,h}$ also exists (and equals $\nu_{\widetilde{\theta_0}}/
\cos (\widetilde{\theta_0} -  \theta) $) and is the limit of $(\phi_n/(nl(A)))_{n\in
  \NN}$ almost surely, even if $\lim_{n\rightarrow \infty} h(n)/n$ does not
exist.

To complete the proof of Corollary \ref{chapitre6corollaire}, it remains to
prove the convergence of $\phi_n/nl(A)$ in $L^1$. Suppose first that
the condition (\ref{conditionO}) is satisfied. Then, one can
find a sequence of sets of edges $(E(n))_{n\in\NN}$ such that for each
$n$, $E(n)$ is a cut between $T(nA, h(n))$ and $B(nA,h(n))$,
$E(n)\subset E(n+1)$ and:
$$\lim_{n\rightarrow \infty}\frac{\card (E(n)) }{n l(A)} \quad
\textrm{exists}\;,$$
cf. Lemma~4.1 in \cite{RossignolTheret08b}, for instance. Now, define:
$$f_n=\frac{\phi_n}{n l(A)}\quad \mbox{ and }\quad g_n= \frac{1}{nl(A)}\sum_{e\in
  E(n)}t(e)\;.$$
We know that $(g_n)_{n\in\NN}$ converges almost surely and in $L^1$, thanks
to the usual law of large numbers, thus the family $(g_n)_{n\in\NN}$ is
equi-integrable. Since $0\leq f_n\leq g_n$ for every $n$, the family
$(f_n)_{n\in\NN}$ is equi-integrable too, so its almost sure convergence
towards $\eta_{\theta,h}$ implies its convergence in $L^1$ towards the same
limit.

\begin{figure}[!ht]
\centering
\begin{picture}(0,0)%
\includegraphics{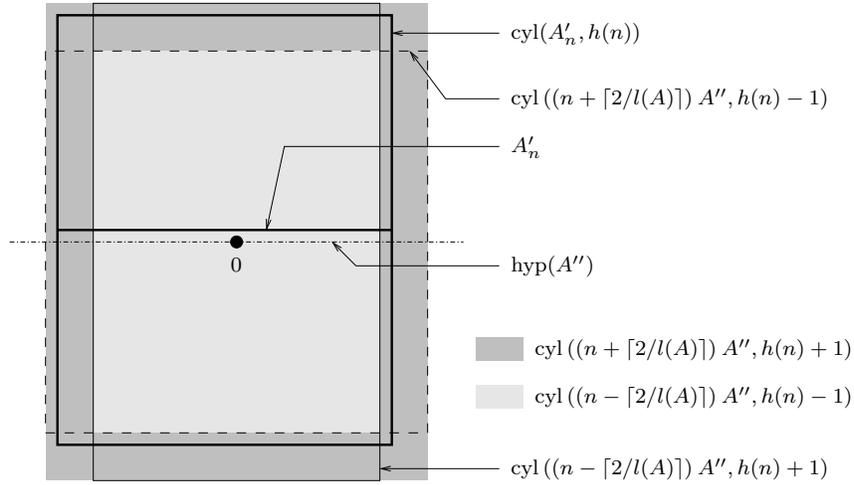}%
\end{picture}%
\setlength{\unitlength}{1973sp}%
\begingroup\makeatletter\ifx\SetFigFont\undefined%
\gdef\SetFigFont#1#2#3#4#5{%
  \reset@font\fontsize{#1}{#2pt}%
  \fontfamily{#3}\fontseries{#4}\fontshape{#5}%
  \selectfont}%
\fi\endgroup%
\begin{picture}(10063,6033)(2539,-6982)
\put(8851,-1411){\makebox(0,0)[lb]{\smash{{\SetFigFont{8}{9.6}{\rmdefault}{\mddefault}{\updefault}{\color[rgb]{0,0,0}$\cyl (A_n',h(n))$}%
}}}}
\put(8851,-2236){\makebox(0,0)[lb]{\smash{{\SetFigFont{8}{9.6}{\rmdefault}{\mddefault}{\updefault}{\color[rgb]{0,0,0}$\cyl \left( \left( n+ \lceil 2/l(A) \rceil \right) A'', h(n)-1 \right)$}%
}}}}
\put(9151,-5386){\makebox(0,0)[lb]{\smash{{\SetFigFont{8}{9.6}{\rmdefault}{\mddefault}{\updefault}{\color[rgb]{0,0,0}$\cyl \left( \left( n+ \lceil 2/l(A) \rceil \right) A'', h(n)+1 \right)$}%
}}}}
\put(9151,-5986){\makebox(0,0)[lb]{\smash{{\SetFigFont{8}{9.6}{\rmdefault}{\mddefault}{\updefault}{\color[rgb]{0,0,0}$\cyl \left( \left( n- \lceil 2/l(A) \rceil \right) A'', h(n)-1 \right)$}%
}}}}
\put(8851,-6886){\makebox(0,0)[lb]{\smash{{\SetFigFont{8}{9.6}{\rmdefault}{\mddefault}{\updefault}{\color[rgb]{0,0,0}$\cyl \left( \left( n- \lceil 2/l(A) \rceil \right) A'', h(n)+1 \right)$}%
}}}}
\put(8851,-2836){\makebox(0,0)[lb]{\smash{{\SetFigFont{8}{9.6}{\rmdefault}{\mddefault}{\updefault}{\color[rgb]{0,0,0}$A_n'$}%
}}}}
\put(8851,-4336){\makebox(0,0)[lb]{\smash{{\SetFigFont{8}{9.6}{\rmdefault}{\mddefault}{\updefault}{\color[rgb]{0,0,0}$\hyp(A'')$}%
}}}}
\put(5401,-4336){\makebox(0,0)[b]{\smash{{\SetFigFont{8}{9.6}{\rmdefault}{\mddefault}{\updefault}{\color[rgb]{0,0,0}$0$}%
}}}}
\end{picture}%
\caption{The cylinder $\cyl(A_n',h(n))$.}
\label{comparaison}
\end{figure}

It remains to show the convergence in $L^1$ without the condition (\ref{conditionO}). Let $A''$ be the translate of $A$ such that $0\in A''$, and
$0$ is the center of $A''$, thus condition (\ref{conditionO}) holds for
$A''$. For any fixed $n$, there exists a segment $A_n'$ which is a translate
of $nA$ by an integer vector and such that $d_\infty(0,nA'_n)<1$ and
$d_\infty(nA'',A'_n)<1$, where $d_\infty$ denotes the distance induced by
$\|.\|_\infty$. We want to compare the maximal flow through
$\cyl(nA'',h(n))$ with the maximal flow through $\cyl (A_n',h(n))$. We have
to distort a little bit the cylinder $\cyl(nA'',h(n))$. We only consider $n$ large enough so that $h(n) >1$. Thus the following inclusions hold:
$$ \cyl \left( \left( n- \left\lceil \frac{2}{l(A)}  \right\rceil
  \right) A'', h(n) -1 \right)\,\subset\, \cyl(A_n',h(n)) \,\subset\, \cyl \left( \left( n+ \left\lceil \frac{2}{l(A)}  \right\rceil
  \right) A'', h(n) +1\right) \,, $$
where $\lceil x \rceil$ is the smallest integer bigger than or equal to
$x$. We get
\begin{align*}
\phi \left( \left( n- \left\lceil \frac{2}{l(A)}  \right\rceil
  \right) A'', h(n) +1\right) \,\leq\,  \phi(A_n',h(n))
\,\leq\, \phi \left(  \left( n+ \left\lceil
      \frac{2}{l(A)}  \right\rceil \right) A'', h(n) -1 \right)\,,
\end{align*}
(see Figure \ref{comparaison}).
Using the convergence in $L^1$ for $A''$ which satisfies the condition
(\ref{conditionO}), we see that
$$\frac{\phi \left( \left( n- \left\lceil \frac{2}{l(A)}  \right\rceil
  \right) A'', h(n) +1\right)}{nl(A)} \quad\textrm{and} \quad
\frac{\phi  \left(  \left( n+ \left\lceil
      \frac{2}{l(A)}  \right\rceil \right) A'', h(n) -1 \right)}{nl(A)}$$ 
converge to
$\eta_{\theta,h}$ in $L^1$ as $n$ goes to infinity. It is obvious that the
small difference in the parameters $n$ and $h(n)$ does not change the value
of the limit $\eta_{\theta,h}$. We get the convergence of
$\tau(A'_n,h(n))/(nl(A))$ to $\eta_{\theta,h}$ in $L^1$. But since
$A_n'$ is an integer translate of $nA$, it implies the convergence of
$\tau(nA,h(n))/(nl(A))$ to $\eta_{\theta,h}$ in $L^1$.

\begin{rem}
\label{remd3}
In dimension $d\geq 3$, if we denote by $\vec{v}$ a unit vector orthogonal
to a non-degenerate hyperrectangle $A$ and by $\overrightarrow{\D_n(A)}$ the set of all
admissible directions for the cylinder $\cyl(nA, h(n))$, i.e., the set
of the vectors $\vec{v}' $ in $S^{d-1}$ such that there exists a
hyperplane $\mathcal{P}$ orthogonal to $\vec{v}' $ that
intersects $\cyl(nA, h(n))$ only on its ``vertical faces'', and if
$\lim_{n\rightarrow \infty} h(n) / n$ exists (thus $\overrightarrow{\D(A)} =
\ad(\underline{\overrightarrow{\D(A)}}) = \ad(\overline{\overrightarrow{\D(A)}})$ exists), we
conjecture that
$$  \lim_{n\rightarrow \infty} \frac{\phi(nA, h(n))}{n^{d-1}\H^{d-1}(A)}
\,=\, \inf_{\vec{v}' \in \overrightarrow{\D(A)}}
\frac{\nu (\vec{v}' )}{|\vec{v} \cdot \vec{v}' |}
\qquad \textrm{a.s.} \,,$$
under assumptions (\ref{chapitre6eq:conditionsFLGN}) on $F$ and if $h(n)$ goes to
infinity with $n$ in such a way that we have $\lim_{n\rightarrow \infty} \log h(n) /
n^{d-1} =0$. We could not prove this conjecture, because we are not able to
prove that $\phi(nA, h(n))$ behaves asymptotically like $\min_{\kappa
  \in K} \phi^{\kappa}(nA, h(n))$, where $K$ is the set of the flat boundary conditions, i.e., the boundary
conditions given by the intersection of a hyperplane with the vertical
faces of $\cyl(nA, h(n))$.
\end{rem}


\def\cprime{$'$}


\end{document}